\documentclass{amsart}
\usepackage{amssymb,latexsym,amsmath,amscd}

\newtheorem{definition}{Definition}[section]

\newtheorem{lem}{Lemma}[section]
\newtheorem{prop}{Proposition}[section] 
\newtheorem{thm}{Theorem}[section]

\newcommand{\pa}{\partial}

\newcommand{\nab}{\nabla}
\newcommand{\ka}{K\"{a}hler }

\newcommand{\BR}{\mathbb{R}}
\newcommand{\BC}{\mathbb{C}}

\newcommand{\BP}{\mathbb{P}}

\newcommand{\CY}{Calabi--Yau }

\newcommand{\we}{\wedge}
\newcommand{\oa}{\omega}
\newcommand{\oal}{\oa^{(1,1)}}
\newcommand{\Oa}{\Omega}
\newcommand{\ga}{\gamma}
\newcommand{\Ga}{\Gamma}
\newcommand{\da}{\delta}
\newcommand{\ta}{\theta}

\newcommand{\ep}{\epsilon}
\newcommand{\al}{\alpha}
\newcommand{\ba}{\beta}
\newcommand{\si}{\sigma}
\newcommand{\Si}{\Sigma}
\newcommand{\Ups}{\Upsilon}
\newcommand{\Lda}{\Lambda}
\newcommand{\lda}{\lambda}
\newcommand{\vi}{\varphi}

\newcommand{\pf}{\emph{Proof. }}
\newcommand{\SM}{S^1 \times M}

\numberwithin{equation}{subsection}

\begin{document}







\begin{center}
\end{center}

\begin{center}
\end{center}

\begin{center} 
\large \bf DESINGULARIZATIONS OF CALABI--YAU 3-FOLDS WITH A CONICAL SINGULARITY \\[0.8cm]
\normalfont{YAT-MING CHAN} \\[0.8cm]
\end{center}

\begin{center}
\bf{Abstract}
\end{center}
\par \small {We study \emph{\CY 3-folds} $M_0$ \emph{with a conical singularity} $x$ modelled on a \emph{\CY cone} $V$. We construct desingularizations of $M_0$, obtaining a 1-parameter family of compact, nonsingular \CY 3-folds which has $M_0$ as the limit. The way we do is to choose an \emph{Asymptotically Conical \CY 3-fold} $Y$ modelled on the same cone $V$, and then glue into $M_0$ at $x$ after applying a homothety to $Y$. We then get a 1-parameter family of \emph{nearly \CY 3-folds} $M_t$ depending on a small real variable $t$. For sufficiently small $t$, we show that the nearly \CY structures on $M_t$ can be deformed to genuine \CY structures, and therefore obtaining the desingularizations of $M_0$. Our result can be applied to resolving orbifold singularities and hence provides a quantitative description of the \CY metrics on the crepant resolutions. \\[0.2cm]

\begin{center}
\large \bf{1. \ Introduction} \\[0.5cm]
\end{center}
\normalsize
\par A \emph{\CY 3-fold} is a \ka manifold \,$(M,J,g)$ of complex dimension 3 with a covariant constant holomorphic volume form $\Oa$ satisfying $\oa^3 = \frac{3i}{4}\,\Oa \wedge \bar{\Oa}$ where $\oa$ is the \ka form of $g$. Equivalently it is a Riemannian 6-fold with a torsion-free SU(3)-structure. \\

\par Our starting point is to introduce the notion of a \emph{nearly \CY structure} on some 6-fold $M$. It is basically an SU(3)-structure with small torsion. We then show that if the torsion is small enough, the nearly \CY structure on $M$ can be deformed to a torsion-free SU(3)-structure, that is, a genuine \CY structure on $M$. This is achieved by applying Joyce's existence result \cite[Thm. 11.6.1]{Joyce1} for torsion-free $G_2$-structures on compact 7-folds. \\

\par We shall then study a kind of singular \CY 3-fold with isolated singularities, known as \emph{\CY 3-folds with conical singularities}. Roughly speaking, they are compact \CY 3-folds with finitely many isolated singularities such that they approach some \emph{\CY cones} near the singular points in some sense. For simplicity we shall consider instead \CY 3-folds with only one conical singularity modelled on a \CY cone, and no other singularities. \\

\par The goal of this paper is to desingularize the compact \CY 3-fold $M_0$ with a conical singularity modelled on a \CY cone $V$ by an analytic technique called \emph{gluing}. The construction proceeds as follows. Suppose $Y$ is an \emph{Asymptotically Conical \CY 3-fold} modelled on the same cone $V$. It has a similar definition to $M_0$, but it is asymptotic to $V$ at infinity. We apply a homothety to $Y$ and then glue into $M_0$ to get a 1-parameter family of compact, nonsingular smooth 6-folds $M_t$. We will then construct \emph{nearly \CY structures} on $M_t$ for small $t > 0$, and when $t$ is sufficiently small, the nearly \CY structures can be deformed to genuine \CY structures by applying our result before. \\

\par An application of our result involves desingularizing \CY 3-orbifolds with isolated singularities. We can then describe what the \CY metrics on the crepant resolution of the orbifold locally look like. \\

\par We begin in \S 2 by giving some background material for this paper. Section 3 introduces nearly \CY structures on 6-dimensional manifolds and the induced $G_2$-structures on 7-dimensional manifolds. We also prove the existence result for genuine \CY structures, using Joyce's result with some modifications. Then we define in $\S$4 the main objects of this paper, namely the \emph{\CY cones, \CY m-folds with a conical singularity} and \emph{Asymptotically Conical \CY m-folds,} and provide some examples. Finally we show in \S5 the construction of the desingularization in the simplest case where there are no obstructions. We shall then give an application of our result and have a short discussion on the obstructed case. \\[0.5cm]   
\textbf{Acknowledgements. } \ I would like to thank my supervisor, Dominic Joyce, for introducing me this topic and for many useful conversations. I would also like to thank The Croucher Foundation for its support. \\[0.3cm]

\begin{center}
\large \bf{2. \ Background material} \\[0.5cm]
\end{center}

\par In this section we provide some background on \CY 3-folds, SU(3)-structures on 6-folds and $G_2$-structures on 7-folds. They will play an essential role in our construction of desingularizations of compact \CY 3-folds with a conical singularity. Let us begin with studying \CY 3-folds. Some useful introductory references on \CY manifolds are \cite{GrossJoyce} and \cite[Chapter 6]{Joyce1}.\\

\renewcommand{\thedefinition}{2.1}
\begin{definition}
\quad \textnormal{A \emph{\CY 3-fold} is a \ka manifold \,$(M,J,g)$ of complex dimension 3 with a covariant constant holomorphic volume form $\Oa$ such that it satisfies $\oa^3 = \frac{3i}{4}\,\Oa \wedge \bar{\Oa}$ where $\oa$ is the \ka form of $g$. We say that $(J,\oa,\Oa)$ constitutes a \emph{\CY structure} on $M$ and write a \CY manifold as a quadruple $(M,J,\oa,\Oa)$.} 
\end{definition}

\par Thus for each $x \in M$, there is an isomorphism between $T_x M$ and $\BC^3$ that identifies $g_x, \oa_x$ and $\Oa_x$ with the flat metric $g_0$, the real 2-form $\oa_0$ and the complex 3-form $\Oa_0$ on $\BC^3$ given  by \\[-0.4cm]
\begin{align*} 
g_0 = |dz_1|^2 + |dz_2|^2 + |dz_3|^2, \quad &\oa_0 = \frac{i}{2}( dz_1 \we d\bar{z}_1 + dz_2 \we d\bar{z}_2 + dz_3 \we d\bar{z}_3) \\
\textnormal{and }\ &\Oa_0 = dz_1 \we dz_2 \we dz_3 , 
\end{align*}
where $(z_1, z_2, z_3)$ are coordinates on $\BC^3$. \CY manifolds are automatically Ricci-flat, and one can use Yau's solution of the Calabi conjecture (\cite{Yau} and \cite[Chapter 5]{Joyce1}) to show the existence of families of \CY manifolds. Another equivalent way of defining a \CY 3-fold is to require that the Riemannian 6-fold $(M,g)$ has holonomy group Hol($g$) contained in SU(3). One can then show that $M$ admits a holomorphic volume form satisfying the normalization formula. \\

\par We shall now consider SU(3)-structures on 6-folds and $G_2$-structures on 7-folds and relate them to \CY structures. An SU(3)-structure on a real 6-fold $M$ is a principal subbundle of the frame bundle of $M$, with fibre SU(3). Each SU(3)-structure gives rise to an almost complex structure $J$, a real 2-form $\oa$ and a complex 3-form $\Oa$ with the properties that 
\begin{enumerate}
	\item $\oa$ is of type (1,1) w.r.t. $J$ and is positive,
	\item $\Oa$ is of type (3,0) w.r.t. $J$ and is nonvanishing, and
	\item $\oa^3 = \frac{3i}{4}\,\Oa \we \bar{\Oa}$. \\[-0.3cm]
\end{enumerate}
We will refer to $(J, \oa, \Oa)$ as an SU(3)-structure. If in addition $d\oa = 0$ and $d\Oa = 0$, then $J$ is integrable and $\Oa$ is a holomorphic (3,0)-form, and the closedness of $\oa$ implies the associated Hermitian metric $g$ is K\"{a}hler. In this case $d\oa$ and $d\Oa$ can be thought of as the \emph{torsion} of the SU(3)-structure, and when they both vanish the SU(3)-structure is \emph{torsion-free}. Note that property (3) implies that the holomorphic (3,0)-form $\Oa$ has constant length, so it is covariant constant. Therefore we have \\
\renewcommand{\theprop}{2.2}
\begin{prop}
\quad Let $M$ be a real 6-fold and $(J, \oa, \Oa)$ an SU(3)-structure on $M$. Let $g$ be the Hermitian metric with Hermitian form $\oa$. Then the followings are equivalent:
\begin{enumerate}
	\item[$(i)$] $d\oa = 0$ \textnormal{and } $d\Oa = 0$ on $M$,
	\item[$(ii)$] $(J, \oa, \Oa)$ is torsion-free,
	\item[$(iii)$] $(J, \oa, \Oa)$ gives a \CY structure on $M$, and
	\item[$(iv)$] \textnormal{Hol($g$) $\subseteq$ SU(3)}. \\
\end{enumerate}
\end{prop}

\par Now we discuss $G_2$-structures in 7-folds. The books by Salamon \cite[$\S$11-$\S$12]{Salamon} and Joyce \cite[Chapter 10]{Joyce1} are good introductions to $G_2$. In the theory of Riemannian holonomy groups, one of the exceptional cases in Berger's classification \cite{Ber} is given by $G_2$ in 7 dimensions. Bryant and Salamon \cite{BS} found explicit, complete metrics with holonomy $G_2$ on noncompact manifolds, and Joyce \cite{Joyce1} constructed examples of compact 7-folds with holonomy $G_2$. The exceptional Lie group $G_2$ is the subgroup of GL(7,$\BR$) preserving the 3-form
\begin{align*} 
\vi_0 &= dx_1 \we dx_2 \we dx_3 + dx_1 \we dx_4 \we dx_5 + dx_1 \we dx_6 \we dx_7 + dx_2 \we dx_4 \we dx_6 \\ 
&\quad - dx_2 \we dx_5 \we dx_7 - dx_3 \we dx_4 \we dx_7 - dx_3 \we dx_5 \we dx_6
\end{align*}
on $\BR^7$ with coordinates $(x_1, \dots , x_7)$. This group also preserves the metric 
$$g_0 = dx_1^2 + \cdots + dx_7^2,$$
the 4-form 
\begin{align*} 
*\vi_0 &= dx_4 \we dx_5 \we dx_6 \we dx_7 + dx_2 \we dx_3 \we dx_6 \we dx_7 + dx_2 \we dx_3 \we dx_4 \we dx_5 \\
&\quad + dx_1 \we dx_3 \we dx_5 \we dx_7 - dx_1 \we dx_3 \we dx_4 \we dx_6 - dx_1 \we dx_2 \we dx_5 \we dx_6 \\
&\quad - dx_1 \we dx_2 \we dx_4 \we dx_7,
\end{align*}
and the orientation on $\BR^7$. Let $X$ be an oriented 7-fold. We say that a 3-form $\vi$ (a 4-form $\psi$) on $X$ is \emph{positive} if for each $p \in X$, there exists an oriented isomorphism between $T_p X$ and $\BR^7$ identifying $\vi$ and the 3-form $\vi_0$ (the 4-form $*\vi_0$). \\

\par A $G_2$-structure on a 7-fold $X$ is a principal subbundle of the oriented frame bundle of $X$, with fibre $G_2$. Thus there is a 1-1 correspondence between positive 3-forms and $G_2$-structures on $X$. Moreover, to any positive 3-form on $X$ one can associate a unique positive 4-form $*\vi$ and metric $g$, such that $\vi$, $*\vi$ and $g$ are identified with $\vi_0$, $*\vi_0$ and $g_0$ under an isomorphism between $T_p X$ and $\BR^7$, for each $p \in M$. We shall refer to $(\vi, g)$ as a $G_2$-structure. Suppose $(\vi, g)$ is a $G_2$-structure on $X$, and $\nab$ is the Levi-Civita connection of $g$. We call $\nab \vi$ the \emph{torsion} of the $G_2$-structure $(\vi, g)$, and if $\nab \vi = 0$, then the $G_2$-structure is \emph{torsion-free}. Here is a result from \cite[Prop. 10.1.3]{Joyce1}:\\

\renewcommand{\theprop}{2.3}
\begin{prop}
\quad Let $X$ be a real 7-fold and $(\vi, g)$ a $G_2$-structure on $X$. Then the followings are equivalent:
\begin{enumerate}
	\item[$(i)$] $(\vi, g)$ is torsion-free,
	\item[$(ii)$] $\nab \vi = 0$ on $X$, where $\nab$ is the Levi-Civita connection of $g$, 
	\item[$(iii)$] $d\vi = d^* \vi = 0$ on $X$, and  
	\item[$(iv)$] \textnormal{Hol($g$)} $\subseteq G_2$, and $\vi$ is the induced 3-form. \\
\end{enumerate}
\end{prop}



\par Now if $(M, J, \oa, \Oa)$ is a \CY 3-fold with the \CY metric $g_M$, then by Proposition 2.2, $(J, \oa, \Oa)$ gives a torsion-free SU(3)-structure on $M$, and Hol($g_M$) $\subseteq$ SU(3). By considering SU(3) as a subgroup of $G_2$, the 7-fold $\SM$ has a torsion-free $G_2$-structure, which is constructed by the following result \cite[Prop. 11.1.2]{Joyce1}: \\

\renewcommand{\theprop}{2.4}
\begin{prop}
\quad Suppose $(J, \oa, \Oa)$ is a torsion-free SU(3)-structure on a 6-fold $M$. Let s be a coordinate on $S^1$. Define a metric $g$ and a 3-form $\vi$ on $\SM$ by
$$ g = ds^2 + g_M \quad and \quad \vi = ds \we \oa + \textnormal{Re}(\Oa). $$
Then $(\vi, g)$ is a torsion-free $G_2$-structure on $\SM$, and 
$$ *\vi = \frac{1}{2}\,\oa \we \oa - ds \we \textnormal{Im}(\Oa). $$ \\
\end{prop}


\begin{center}
\large \bf{3. \ Nearly Calabi--Yau structures} \\[0.5cm]
\end{center}

\normalsize

\par This section introduces the notion of a \emph{nearly \CY structure} on an oriented 6-fold $M$. We begin in $\S$3.1 by giving the definition of a nearly \CY structure $(\oa, \Oa)$ on $M$, and showing that if $M$ admits a genuine \CY structure, then any $(\oa, \Oa)$ which is sufficiently close to it gives a nearly \CY structure. Section 3.2 constructs $G_2$-structures on the 7-fold $\SM$. Finally, we give the main result of the section, the existence of genuine \CY structures on $M$, in $\S$3.3. It is based on the existence result for torsion-free $G_2$-structures on compact 7-folds by Joyce \cite[Thm. 11.6.1]{Joyce1}. \\

\begin{center}
\bf{3.1. \ Introduction to nearly \CY structures}\\[0.3cm]
\end{center}

\par Let $M$ be an oriented 6-fold. A \emph{nearly \CY structure} on $M$ consists of a real closed 2-form $\oa$, and a complex closed 3-form $\Oa$ on $M$. Basically, the idea of a nearly \CY structure $(\oa, \Oa)$ is that it corresponds to an SU(3)-structure with ``small torsion", and hence approximates a genuine \CY structure, which is equivalent to a torsion-free SU(3)-structure. So let us start with generating an SU(3)-structure on $M$ from $(\oa, \Oa)$. \\

\par First we write $\Oa = \ta_1 + i \ta_2$, so $\ta_1$ and $\ta_2$ are both real closed 3-forms. Suppose $\ta_1$ has stabilizer SL(3,\,$\BC$) $\subset$ GL$_+$(6,\,$\BR$) at each  $p \in M$, then the orbit of $\ta_1$ in $\bigwedge^3 T^{*}_{p}M$ under the action of GL$_+$(6,\,$\BR$) is GL$_+$(6,\,$\BR$)/SL(3,\,$\BC$). For each $p \in M$, define $\bigwedge^{3}_{+} T^{*}_{p}M$ to be the subset of 3-forms $\ta \in \bigwedge^{3} T^{*}_{p}M$ for which there exists an oriented isomorphism between $T_p M$ and $\BR^6 \cong \BC^3$ identifying $\ta$ and the 3-form Re$(dz_1 \we dz_2 \we dz_3)$ where $(z_1, z_2, z_3)$ are coordinates on $\BC^3$. Then $\bigwedge^{3}_{+} T^{*}_{p}M \cong$ GL$_+$(6,\,$\BR$)/SL(3,\,$\BC$), as Re$(dz_1 \we dz_2 \we dz_3)$ has stabilizer SL(3,\,$\BC$). Then $\ta_1 |_p$ lies in $\bigwedge^{3}_{+} T^{*}_{p}M$ for each $p \in M$. It is easy to see that dim $\bigwedge^{3}_{+} T^{*}_{p}M$ = dim GL$_+$(6,\,$\BR$)/SL(3,\,$\BC$) = dim $\bigwedge^3 T^{*}_{p}M$ = 20, so $\bigwedge^{3}_{+} T^{*}_{p}M$ is an open subset of $\bigwedge^3 T^{*}_{p}M$ for each $p \in M$. Therefore any 3-form on $M$ which is sufficiently close to a 3-form in $\bigwedge^{3}_{+} T^{*}_{p}M$ still lies in $\bigwedge^{3}_{+} T^{*}_{p}M$, or equivalently, has stabilizer SL(3,\,$\BC$) at each point on $M$.  \\

\par The oriented frame bundle $F_+$ of $M$ is the bundle over $M$ whose fibre at $p \in M$ is the set of oriented isomorphisms between $T_p M$ and $\BR^6$. Let $P$ be the subset of $F_+$ consisting of oriented isomorphisms between $T_p M$ and $\BR^6$ which identify $\ta_1$ at $p$ and Re$(dz_1 \we dz_2 \we dz_3)$. It is well-defined as we have assumed that $\ta_1 |_p \in \bigwedge^{3}_{+} T^{*}_{p}M$. Thus $\ta_1$ defines a principal subbundle $P$ of the oriented frame bundle $F_+$, with fibre SL(3,\,$\BC$), that is, an SL(3,\,$\BC$)-structure on $M$. As SL(3,\,$\BC$) acts on $\BR^6 \cong \BC^3$ preserving the complex structure $J_0$ on $\BC^3$, we obtain a unique almost complex structure $J'$ on $M$. \\

\par Note that the forms $\oa$, $\Oa$ are not necessarily of type (1,1) and (3,0) with respect to $J'$ respectively. We then denote $\oal$ by the (1,1)-component of $\oa$ with respect to $J'$ and define a 3-form $\ta_2 '$ on $M$ by $\ta_2'(u,v,w) := \ta_1 (J'\,u,v,w)$ for all $u,v,w \in TM$, or in index notation, $ (\ta_2')_{abc} = (J')^{d}_{a}\, (\ta_1)_{dbc}$. Suppose that $\oal$ is a positive (1,1)-form, that is, $\oal(v, J'v) > 0$ for any nonzero $v \in TM$. Write $\Oa' = \ta_1 + i \ta_2'$, then $\Oa'$ is a (3,0)-form with respect to $J'$. In general, $\ta_2'$ will not be a closed 3-form, unless $J'$ is integrable. \\

\par We want $(J', \oal, \Oa')$ to be an SU(3)-structure, but the problem with this is the usual normalization formula defining a \CY manifold may not hold for $\oal$ and $\Oa'$, that is, $(\oal)^3 \neq \frac{3i}{4}\,\Oa' \we \bar{\Oa}'\, (= \frac{3}{2}\,\ta_1 \we \ta_2')$ in general. We then define a smooth function $f : M \rightarrow (0, \infty)$ by
\renewcommand{\theequation}{3.1}
	\begin{equation}
	(\oal)^3 = f \cdot \frac{3}{2}\,\ta_1 \we \ta_2'. \\[-0.2cm] 
	\end{equation}
Consequently, if we rescale $\oal$ by setting $\oa' = f^{-\frac{1}{3}}\,\oal $, we have $\oa'^3 = \frac{3}{2}\,\ta_1 \we \ta_2'$. Given that $\oal$, and hence $\oa'$, is positive, then one can determine a Hermitian metric $g_M$ on $M$ from $\oa'$ and $J'$ by $g_M (u,v) = \oa' (u, J'v)$ for all $u,v \in TM$. \\

\par Now we are ready to give the definition of a nearly \CY structure on $M$: \\

\renewcommand{\thedefinition}{3.1}
\begin{definition}
\quad \textnormal{Let $M$ be an oriented 6-fold, and let $\oa$ be a real closed 2-form, and $\Oa = \ta_1 + i\ta_2$ a complex closed 3-form on $M$. Let $\ep_0 \in (0,1]$ be a fixed small constant, to be chosen later in Lemma 3.3. Then $(\oa, \Oa)$ constitutes a \emph{nearly \CY structure} on $M$ if}

\begin{itemize}
	\item[\textnormal{(i)}] \textnormal{the real closed 3-form $\ta_1$ has stabilizer SL(3,\,$\BC$) at each point of $M$, or equivalently, it lies in $\bigwedge^{3}_{+} T^{*}_{p}M$ for each $p \in M$. \\[0.05cm] Then we can associate a unique almost complex structure $J'$ and a unique real 3-form $\ta_2'$ such that $\Oa' = \ta_1 + i\ta_2'$ is a (3,0)-form with respect to $J'$. \\[-0.3cm]}
	
	\item[\textnormal{(ii)}] \textnormal{the (1,1)-component $\oal$ of $\oa$ with respect to $J'$ is positive. \\[0.05cm] Then we can associate a Hermitian metric $g_M$ on $M$ from $\oa'$ and $J'$, where $\oa' = f^{-\frac{1}{3}}\,\oal $ is the rescaled (1,1)-part of $\oa$ and $f$ is defined by (3.1). \\[-0.3cm]}
	
	\item[\textnormal{(iii)}] \textnormal{the following inequalities hold for some $\ep$ with $0< \ep \leq \ep_0$:}
	\renewcommand{\theequation}{3.2} 
	\begin{align} 	\left| \ta_2 - \ta_2' \right|_{g_M} &< \ep	\, ,\\ \tag{3.3}
	| \oa^{(2,0)} |_{g_M} &< \ep	\, , \textnormal{and} \\ \tag{3.4}
	\left| \oa^3 - \frac{3}{2}\, \ta_1 \we \ta_2 \right|_{g_M} &< \ep  \notag
	\end{align} 
\textnormal{where the norms $|\cdot|_{g_M}$ are measured by the metric $g_M$. \\[-0.2cm]}	
\end{itemize}
\end{definition} 

\par If $(\oa, \Oa)$ is a nearly \CY structure on $M$, one can show that the function $f$ defined in (3.1) satisfies 
\renewcommand{\theequation}{3.5}
	\begin{align}
	| f - 1 | < C_0 \ep, 
	\end{align}
for some constant $C_0 > 0$, i.e. $f$ is approximately equal to 1 for sufficiently small $\ep$, as we would expect. \\

\par The next result shows that if $M$ admits a genuine \CY structure, then any $(\oa, \Oa)$ which is sufficiently close to it gives a nearly \CY structure on $M$. \\

\renewcommand{\theprop}{3.2} 
\begin{prop}
\quad There exist constants $\ep_1, C, C' >0$ such that whenever $0 < \ep \leq \ep_1$, the following is true. \\[-0.3cm]
\par Let $M$ be an oriented 6-fold. Suppose $(\tilde{J}, \tilde{\oa}, \tilde{\Oa})$ is a \CY structure with \CY metric $\tilde{g}$, $\oa$ a real closed 2-form, and $\Oa = \ta_1 + i \ta_2$ a complex closed 3-form on $M$, satisfying \renewcommand{\theequation}{3.6} 	
	\begin{align}
	| \tilde{\oa} - \oa |_{\tilde{g}} < \ep \quad and \quad | \tilde{\Oa} - \Oa |_{\tilde{g}} < \ep,
	\end{align}
then $(\oa, \Oa)$ is a nearly \CY structure on $M$ with metric $g_M$ satisfying 
\renewcommand{\theequation}{3.7}
	\begin{align}
	| \tilde{g} - g_M |_{\tilde{g}} < C \ep \quad and \quad | \tilde{g}^{-1} - g^{-1}_{M} |_{\tilde{g}} < C' \ep. 
	\end{align}
\end{prop}
\pf \ \,From (3.6) we have $|\textnormal{Re}(\tilde{\Oa}) - \ta_1 |_{\tilde{g}} < \ep $, which means that if we choose $\ep_1$ to be sufficiently small, then $\ta_1$ has stabilizer SL(3,$\BC$) since the stabilizer condition is an open condition as we mentioned before. So we can associate a unique almost complex structure $J'$, with $| \tilde{J} - J' |_{\tilde{g}} < C_1 \ep$ for some constant $C_1 > 0$, and a unique real 3-form $\ta_2'$ such that $\Oa' = \ta_1 + i\ta_2'$ is a (3,0)-form with respect to $J'$. \\

\par One can deduce from (3.6) and $| \tilde{J} - J' |_{\tilde{g}} < C_1 \ep$ that $| \tilde{\oa} - \oa^{(1,1)} |_{\tilde{g}} < C_2 \ep$ for some $C_2 > 0$. Make $\ep_1$ smaller if necessary, then $\oa^{(1,1)}$ is a positive (1,1)-form with respect to $J'$ since the positivity is also an open condition. Then we can define a metric $g_M$ by $g_M (u,v) = \oa' (u, J'v)$ for all $u,v \in TM$, where $\oa' = f^{-\frac{1}{3}}\,\oal $ and $f$ is defined by (3.1). Now we show that $f$ is close to 1. In fact, 
\begin{align} 
| (f-1) \oa'^3 |_{\tilde{g}} &= | (\oa^{(1,1)})^3 - \oa'^3 |_{\tilde{g}} \notag \\ \notag
&\leq | (\oa^{(1,1)})^3 - \tilde{\oa}^3 |_{\tilde{g}} + \left| \tilde{\oa}^3 - \frac{3}{2} \textnormal{Re}(\tilde{\Oa}) \we \textnormal{Im}(\tilde{\Oa}) \right|_{\tilde{g}} \\ \notag
&\quad + \left| \frac{3}{2} \textnormal{Re}(\tilde{\Oa}) \we \textnormal{Im}(\tilde{\Oa}) - \frac{3}{2} \ta_1 \we \ta_2' \right|_{\tilde{g}}. \tag{3.8} \\[-0.3cm] \notag
\end{align}
Since $| \tilde{\oa} - \oa^{(1,1)} |_{\tilde{g}} < C_2 \ep$, the first term of right hand side of (3.8) is of size $O(\ep)$. The second term vanishes as $(\tilde{\oa}, \tilde{\Oa})$ is a \CY structure. From $|\textnormal{Re}(\tilde{\Oa}) - \ta_1 |_{\tilde{g}} < \ep $, and $| \tilde{J} - J' |_{\tilde{g}} < C_1 \ep$, we have  
$|\textnormal{Im}(\tilde{\Oa}) - \ta_2' |_{\tilde{g}} < C_3 \ep $ for some $C_3 > 0$ and hence the third term also has size $O(\ep)$. Summing up, we have $f-1$ is of size $O(\ep)$. \\




\par It can then be shown that $| \tilde{\oa} - \oa' |_{\tilde{g}} < C_4 \ep$ for some $C_4 > 0$. Together with $| \tilde{J} - J' |_{\tilde{g}} < C_1 \ep$, we obtain first part of (3.7), that is, $| \tilde{g} - g_M |_{\tilde{g}} < C \ep$ for some $C > 0$. If $\ep$ is small enough such that $C \ep < \frac{1}{2}$, then one can deduce that $| \tilde{g}^{-1} - g^{-1}_{M} |_{\tilde{g}} < C' \ep$ for some $C' > 0$. This implies that $\tilde{g}$ and $g_M$ are uniformly equivalent metrics, and hence norms of any tensor on $M$ taken with respect to $\tilde{g}$ and with respect to $g_M$ differ by a bounded factor. \\

\par It remains to check (3.2)-(3.4) of Definition 3.1. But it is not hard to get bounds for (3.2)-(3.4) in terms of $\ep$ with respect to the metric $\tilde{g}$, and so by making $\ep$ smaller and using the equivalence of the metrics, we obtain (3.2)-(3.4). \hfill $\Box$ \\[0.7cm]

\begin{center}
\bf{3.2. \ $G_2$-structures on $\SM$} \\[0.3cm]
\end{center}

\par Let $(\oa, \Oa)$ be a nearly \CY structure on $M$. From $\S$3.1 we know that $(J', \oa', \Oa')$ gives an SU(3)-structure with metric $g_M$ on $M$. In this section, we would like to discuss $G_2$-structures on the 7-fold $\SM$, which is essential for the main result in next section. Let $s$ be a coordinate on $S^1$. Now define a 3-form $\vi'$ and a metric $g'$ on $\SM$ by
\renewcommand{\theequation}{3.9}
	\begin{align} 
	\vi' = ds \we \oa' + \ta_1 \quad \textnormal{and} \quad g' = ds^2 + g_M. 
	\end{align}
It turns out that $(\vi', g')$ defines a $G_2$-structure (with torsion) on $\SM$. The associated 4-form $\ast_{g'} \vi'$ on $\SM$ is then given by \renewcommand{\theequation}{3.10}
	\begin{align} 
	\ast_{g'} \vi' = \frac{1}{2}\,\oa' \we \oa' - ds \we \ta_2'.
	\end{align}
Also, we can construct another 3-form $\vi$ and 4-form $\chi$ on $\SM$ by \renewcommand{\theequation}{3.11}
	\begin{align} 
	\vi = ds \we \oa + \ta_1 \quad \textnormal{and} \quad \chi = \frac{1}{2}\,\oa \we \oa - ds \we \ta_2.
	\end{align}
The next lemma shows that the forms in (3.11) are close to the $G_2$-forms $\vi'$ and $\ast_{g'} \vi'$ if we make $\ep_0$ in the definition of nearly \CY manifolds to be sufficiently small. \\

\renewcommand{\thelem}{3.3} 
\begin{lem}
\quad There exist constants $C_1, C_2, C_3$ and $C_4 > 0$ such that if $\ep_0$ in Definition 3.1 is chosen sufficiently small, then the following is true. \\[-0.3cm]

\par Let $(\vi', g')$ be the $G_2$-structure given by (3.9), $\ast_{g'} \vi'$ the associated 4-form given by (3.10), $\vi$ the 3-form and $\chi$ the 4-form given by (3.11) on $\SM$. Then \renewcommand{\theequation}{3.12}      			\begin{align}	
	|\vi - \vi'|_{g'} < C_1 \ep \\[-0.5cm] \notag
	\end{align}
where $\ep \in (0, \ep_0]$ is the small constant in (iii) of Definition 3.1. Hence $\vi$ is also a positive 3-form on $\SM$, and it defines another $G_2$-structure $(\vi, g)$. Moreover, the associated metric $g$ and the 4-form $\ast_{g}\vi$ satisfy \renewcommand{\theequation}{3.13}
	\begin{align} 
	|g - g' |_{g'} < C_2 \ep, \quad |g^{-1} - g'^{-1} |_{g'} &< C_3 \ep \ and \\ \tag{3.14}
	|\!\ast_{g}\!\vi - \chi |_{g'} < C_4 \ep. \\ \notag
	\end{align} 
\end{lem}
\pf \ \,From (3.9) and (3.11), we have $| \vi - \vi' |_{g'} = | ds \we (\oa - \oa') |_{g'} = | \oa - \oa' |_{g_M} < C_1\ep$ for some constant $C_1 > 0$, where we used (3.3) and (3.5). Now we choose $\ep_0$ in Definition 3.1 such that $C_1 \ep_0$ is small enough and $\vi$ is a positive 3-form on $\SM$. \\

\par It follows from general facts about $G_2$-forms that the associated metric $g$ satisfies $|g - g' |_{g'} < C_2 \ep$ for some $C_2 > 0$. Also, by using the same argument as in Proposition 3.2, we obtain $|g^{-1} - g'^{-1} |_{g'} < C_3 \ep$ for some $C_3 > 0$. Now, \renewcommand{\theequation}{3.15}
\begin{align}
|\!\ast_{g}\!\vi - \chi |_{g'} &\leq |\!\ast_{g}\!\vi - \ast_{g'} \vi' |_{g'} + |\!\ast_{g'}\!\vi' - \chi |_{g'} \notag \\ 
&\leq |(\ast_{g} - \ast_{g'}) \vi |_{g'} + |\!\ast_{g'} (\vi - \vi') |_{g'} + |\!\ast_{g'}\!\vi' - \chi |_{g'}.
\end{align}
The first term of right hand side has the same size as $|g - g' |_{g'}$, which is bounded in (3.13). The second term is just $|\vi - \vi'|_{g'}$ since $\ast_{g'}$ is an isometry with respect to $g'$, and so it is bounded in (3.12). For the last term, we have 
$$ |\!\ast_{g'}\!\vi' - \chi |_{g'} = \left| \frac{1}{2}\,(\oa' \we \oa' - \oa \we \oa ) - ds \we ( \ta_2' - \ta_2 ) \right|_{g'} $$
from (3.10) and (3.11), which can be shown by using $| \oa - \oa' |_{g_M} < C_1 \ep$ and (3.2) that has size $O(\ep)$. Summing up together, we get $|\!\ast_{g}\!\vi - \chi |_{g'} < C_4 \ep$ for some constant $C_4 > 0$. \hfill $\Box$ \\[0.5cm] 


\begin{center}
\bf{3.3. \ An existence result for \CY structures on $M$} \\[0.3cm]
\end{center}

\par In the last part of this section we present our main result which shows that when $\ep_0$ is sufficiently small, the nearly \CY structure on $M$ can be deformed to a genuine \CY structure. The proof of it is based on an existence result for torsion-free $G_2$-structures given by Joyce \cite[Thm. 11.6.1]{Joyce1}, which shows using analysis that any $G_2$-structure on a compact 7-fold with sufficiently small torsion can be deformed to a nearby torsion-free $G_2$-structure. We shall adopt a slightly modified version of this result, which improves the bounds of various norms to fit into our situation. We refer to the spaces $L^q$, $L^{q}_{k}$, $C^{k}$ and $C^{k, \al}$ as the Banach spaces defined in \cite[$\S$1.2]{Joyce1}. We begin by stating Joyce's result, with improvements to powers of $t$: \\

\renewcommand{\thethm}{3.4} 
\begin{thm}
\quad Let $\kappa >0$ and $D_1, D_2, D_3 > 0$ be constants. Then there exist constants $\ep \in (0,1]$ and $K > 0$ such that whenever $0 < t \leq \ep$, the following is true. \\[-0.4cm]
\par Let $X$ be compact 7-fold, and $(\vi, g)$ a $G_2$-structure on $X$ with $d\vi = 0$. Suppose $\psi$ is a smooth 3-form on $X$ with $d^* \psi = d^* \vi$, and 
\begin{itemize}
	\item[\textnormal{(i)}] $\| \psi \|_{L^2} \leq D_1 t^{\frac{7}{2}+\kappa}$, $\| \psi \|_{C^0} \leq D_1 t^{\kappa}$\ \,and\ \,$\| d^* \psi \|_{L^{14}} \leq D_1 t^{-\frac{1}{2}+\kappa}$,
	\item[\textnormal{(ii)}] the injectivity radius $\da(g)$ satisfies $\da(g) \geq D_2 t$, and
	\item[\textnormal{(iii)}] the Riemann curvature $R(g)$ satisfies $\| R(g) \|_{C^0} \leq D_3 t^{-2}$. \\[-0.5cm]
\end{itemize}
Then there exists a smooth, torsion-free $G_2$-structure $(\tilde{\vi}, \tilde{g})$ on $X$ such that $\| \tilde{\vi} - \vi \|_{C^0} \leq K t^{\kappa}$\, and\, $[\tilde{\vi}] = [\vi]$ in $H^3 (X, \BR)$. \\ 
\end{thm}

\par The proof of it depends upon the following two results. We state them here and then we will improve the powers of $t$ so that Theorem 3.4 can be modified to fit into our situation for the 7-fold $\SM$. \\

\renewcommand{\thethm}{3.5} 
\begin{thm}
\quad Let $D_2, D_3, t > 0$ be constants, and suppose $(X,g)$ is a complete Riemannian 7-fold, whose injectivity radius $\da(g)$ and Riemann curvature $R(g)$ satisfy $\da(g) \geq D_2 t$ and $\| R(g) \|_{C^0} \leq D_3 t^{-2}$. Then there exist $K_1, K_2 > 0$ depending only on $D_2$ and $D_3$, such that if \,$\chi \in L^{14}_{1} (\Lambda^3 T^* X) \cap L^2 (\Lambda^3 T^* X)$ then 
\begin{align*}
\| \nab \chi \|_{L^{14}} &\leq K_1 \,(\| d\chi \|_{L^{14}} + \| d^* \chi \|_{L^{14}} + t^{-4} \| \chi \|_{L^2} ) \\[0.2cm]
and \quad \| \chi \|_{C^0} &\leq K_2 \,(t^{\frac{1}{2}} \| \nab \chi \|_{L^{14}} + t^{-\frac{7}{2}} \| \chi \|_{L^2}).
\end{align*}
\end{thm}

\par The second result is: \\

\renewcommand{\thethm}{3.6} 
\begin{thm}
\quad Let $\kappa > 0$ and $D_1, K_1, K_2 > 0$ be constants. Then there exist constants $\ep \in (0,1]$, $K_3$ and $K > 0$ such that whenever $0 < t \leq \ep$, the following is true. \\[-0.4cm]
\par Let $X$ be a compact 7-fold, and $(\vi, g)$ a $G_2$-structure on $X$ with $d\vi = 0$. Suppose $\psi$ is a smooth 3-form on $X$ with $d^* \psi = d^* \vi$, and 
\begin{itemize}
	\item[\textnormal{(i)}] $\| \psi \|_{L^2} \leq D_1 t^{\frac{7}{2}+\kappa}$, $\| \psi \|_{C^0} \leq D_1 t^{\kappa}$\ \,and\ \,$\| d^* \psi \|_{L^{14}} \leq D_1 t^{-\frac{1}{2}+\kappa}$,
	\item[\textnormal{(ii)}] if $\chi \in L^{14}_{1} (\Lambda^3 T^* X)$ then $\| \nab \chi \|_{L^{14}} \leq K_1 \,(\| d\chi \|_{L^{14}} + \| d^* \chi \|_{L^{14}} + t^{-4} \| \chi \|_{L^2} )$, 
	\item[\textnormal{(iii)}] if $\chi \in L^{14}_{1} (\Lambda^3 T^*X)$ then $\| \chi \|_{C^0} \leq K_2 \,(t^{\frac{1}{2}} \| \nab \chi \|_{L^{14}} + t^{-\frac{7}{2}} \| \chi \|_{L^2})$.
	\end{itemize}
Let $\ep_1$ be as in Definition 10.3.3, and $F$ as in Proposition 10.3.5 of Joyce \cite{Joyce1}. Denote $\pi_1$ by the orthogonal projection from $\Lda^3 T^* X$ to the 1-dimensional component of the irreducible representation of $G_2$. Then there exist sequences $\{\eta_j\}^{\infty}_{j=0}$ in $L^{14}_{2} (\Lambda^2 T^*X)$ and $\{f_j\}^{\infty}_{j=0}$ in $L^{14}_{1} (X)$ with $\eta_0 = f_0 = 0$, satisfying the equations 
$$ (dd^* + d^* d)\eta_j = d^* \psi + d^* (f_{j-1} \psi) + \ast \, dF(d\eta_{j-1}) \ and \ f_j \vi = \frac{7}{3} \pi_1 (d\eta_j) $$ 
for each $j >0$, and the inequalities
\begin{align*}
&(a)\ \| d\eta_j \|_{L^2} \leq 2 D_1 t^{\frac{7}{2} + \kappa},  &&(d)\ \| d\eta_j - d\eta_{j-1} \|_{L^2} \leq 2 D_1 2^{-j} t^{\frac{7}{2} + \kappa},  \\ 
&(b)\ \| \nab d\eta_j \|_{L^{14}} \leq K_3 t^{-\frac{1}{2} + \kappa}, &&(e)\ \| \nab (d\eta_j - d\eta_{j-1}) \|_{L^{14}} \leq K_3 2^{-j} t^{-\frac{1}{2} + \kappa}, \\
&(c)\ \| d\eta_j \|_{C^0} \leq K t^{\kappa} \leq \ep_1 \qquad and \ &&(f)\ \| d\eta_j - d\eta_{j-1} \|_{C^0} \leq K 2^{-j} t^{\kappa}. \\
\end{align*}
\end{thm}

\par We shall first modify Theorem 3.5 by considering the 6-dimensional version of those analytic estimates. In Theorem 3.5, the first inequality is derived from an elliptic regularity estimate for the operator $d + d^*$ on 3-forms on $X$. The second inequality follows from the Sobolev Embedding Theorem, which states that $L^{q}_{k}$ embeds in $C^{l,\al}$ if $\frac{1}{q} \leq \frac{k-l-\al}{n}$ where $n$ is the dimension of the underlying Riemannian manifold. For the 7-dimensional case, we have $L^{14}_{1}$ embeds in $C^{0,1/2}$ which then embeds in $C^{0}$, whereas in 6 dimensions, we have $L^{12}_{1}$ embeds in $C^{0,1/2}$. We can use this to show \\

\renewcommand{\thethm}{3.7} 
\begin{thm}
\quad Let $D_2, D_3, t > 0$ be constants, and suppose $(M,g)$ is a complete Riemannian 6-fold, whose injectivity radius $\da(g)$ and Riemann curvature $R(g)$ satisfy $\da(g) \geq D_2 t$ and $\| R(g) \|_{C^0} \leq D_3 t^{-2}$. Then there exist $K_1, K_2 > 0$ depending only on $D_2$ and $D_3$, such that if \,$\chi \in L^{12}_{1} (\Lda^3 T^* M) \cap L^2 (\Lda^3 T^* M)$ then 
\begin{align*}
\| \nab \chi \|_{L^{12}} &\leq K_1 \,(\| d\chi \|_{L^{12}} + \| d^* \chi \|_{L^{12}} + t^{-\frac{7}{2}} \| \chi \|_{L^2} ) \\[0.2cm]
and \quad \| \chi \|_{C^0} &\leq K_2 \,(t^{\frac{1}{2}} \| \nab \chi \|_{L^{12}} + t^{-3} \| \chi \|_{L^2}). \\
\end{align*} 
\end{thm}  

\par The proof of it is similar to \cite[Thm. G1, p.298]{Joyce1}. We can first prove the case for $t=1$, and the case for general $t > 0$ follows by conformal rescaling: apply the $t = 1$ case to the metric $t^{-2}g$. The factors of $t$ compensate for powers of $t$ which the norms scaled by in replacing $g$ by $t^{-2}g$. \\ 


\par By considering $S^1$-invariant forms and $S^1$-invariant $G_2$-structures on the 7-fold $\SM$, we can use the Sobolev Embedding Theorem in 6 dimensions, rather than in 7 dimensions. This has an advantage that the powers of $t$ and the inequalities are calculated in 6 dimensions, though the norms are computed on the 7-fold $\SM$. Here is the modified version of Theorem 3.6: \\

\renewcommand{\thethm}{3.8} 
\begin{thm}
\quad Let $\kappa > 0$ and $D_1, K_1, K_2 > 0$ be constants. Then there exist constants $\ep \in (0,1]$, $K_3$ and $K > 0$ such that whenever $0 < t \leq \ep$, the following is true. \\[-0.4cm]
\par Let $M$ be a compact 6-fold, and $(\vi, g)$ an $S^1$-invariant $G_2$-structure on $\SM$ with $d\vi = 0$. Suppose $\psi$ is an $S^1$-invariant smooth 3-form on the 7-fold $\SM$ with $d^* \psi = d^* \vi$, and
\begin{itemize}
	\item[\textnormal{(i)}] $\| \psi \|_{L^2} \leq D_1 t^{3+\kappa}$, $\| \psi \|_{C^0} \leq D_1 t^{\kappa}$\ \,and\ \,$\| d^* \psi \|_{L^{12}} \leq D_1 t^{-\frac{1}{2}+\kappa}$,
	\item[\textnormal{(ii)}] if $\chi \in L^{12}_{1} (\Lda^3 T^* (\SM))$ is $S^1$-invariant, then $\| \nab \chi \|_{L^{12}} \leq K_1 \,(\| d\chi \|_{L^{12}} + \| d^* \chi \|_{L^{12}} + t^{-\frac{7}{2}} \| \chi \|_{L^2} )$, 
	\item[\textnormal{(iii)}] if $\chi \in L^{12}_{1} (\Lda^3 T^* (\SM))$ is $S^1$-invariant, then $\| \chi \|_{C^0} \leq K_2 \,(t^{\frac{1}{2}} \| \nab \chi \|_{L^{12}} + t^{-3} \| \chi \|_{L^2})$.
	\end{itemize} 
With the same notations as in Theorem 3.6, there exist sequences $\{\eta_j\}^{\infty}_{j=0}$ in $L^{12}_{2} (\Lda^2 T^*(\SM))$ and $\{f_j\}^{\infty}_{j=0}$ in $L^{12}_{1} (\SM)$ with $\eta_j, f_j$ being all $S^1$-invariant and $\eta_0 = f_0 = 0$, satisfying the equations 
$$ (dd^* + d^* d)\eta_j = d^* \psi + d^* (f_{j-1} \psi) + \ast \, dF(d\eta_{j-1}) \ and \ f_j \vi = \frac{7}{3} \pi_1 (d\eta_j) $$ 
for each $j >0$, and the inequalities
\begin{align*}
&(a)\ \| d\eta_j \|_{L^2} \leq 2 D_1 t^{3 + \kappa},  &&(d)\ \| d\eta_j - d\eta_{j-1} \|_{L^2} \leq 2 D_1 2^{-j} t^{3 + \kappa},  \\ 
&(b)\ \| \nab d\eta_j \|_{L^{12}} \leq K_3 t^{-\frac{1}{2} + \kappa}, &&(e)\ \| \nab (d\eta_j - d\eta_{j-1}) \|_{L^{12}} \leq K_3 2^{-j} t^{-\frac{1}{2} + \kappa}, \\
&(c)\ \| d\eta_j \|_{C^0} \leq K t^{\kappa} \leq \ep_1 \qquad and \ &&(f)\ \| d\eta_j - d\eta_{j-1} \|_{C^0} \leq K 2^{-j} t^{\kappa}. 
\end{align*}
Here $\nab$ and $\| \cdot \|$ are computed using $g$ on $\SM$. \\
\end{thm}  

\par Thus Theorem 3.8 is essentially an $S^1$-invariant version of Theorem 3.6. We are now ready to state the modified version of Theorem 3.4, to be used to obtain our main result. \\

\renewcommand{\thethm}{3.9} 
\begin{thm}
\quad Let $\kappa >0$ and $D_1, D_2, D_3 > 0$ be constants. Then there exist constants $\ep \in (0,1]$ and $K > 0$ such that whenever $0 < t \leq \ep$, the following is true. \\[-0.4cm]
\par Let $M$ be a compact 6-fold, and $(\vi, g)$ an $S^1$-invariant $G_2$-structure on $\SM$ with $d\vi = 0$. Suppose $\psi$ is an $S^1$-invariant smooth 3-form on $\SM$ with $d^* \psi = d^* \vi$, and
\begin{itemize}
	\item[\textnormal{(i)}] $\| \psi \|_{L^2} \leq D_1 t^{3+\kappa}$, $\| \psi \|_{C^0} \leq D_1 t^{\kappa}$\ \,and\ \,$\| d^* \psi \|_{L^{12}} \leq D_1 t^{-\frac{1}{2}+\kappa}$,
	\item[\textnormal{(ii)}] the injectivity radius $\da(g)$ satisfies $\da(g) \geq D_2 t$, and
	\item[\textnormal{(iii)}] the Riemann curvature $R(g)$ satisfies $\| R(g) \|_{C^0} \leq D_3 t^{-2}$. \\[-0.5cm]
\end{itemize}
Then there exists a smooth, torsion-free $G_2$-structure $(\tilde{\vi}, \tilde{g})$ on $\SM$ such that $\| \tilde{\vi} - \vi \|_{C^0} \leq K t^{\kappa}$\, and\, $[\tilde{\vi}] = [\vi]$ in $H^3 (\SM, \BR)$. \\ 
\end{thm} 

\par Theorem 3.9 follows from Theorems 3.7 and 3.8 as on \cite[p.296-297]{Joyce1}. \\

\par In the remaining part of this section, we shall derive an existence result for genuine \CY structures. Our strategy is the following. We start with the nearly \CY structure $(\oa, \Oa)$ on $M$, then from $\S$3.2 one can induce a $G_2$-structure $(\vi, g)$ on $\SM$. It can then be shown in the following theorem that, under appropriate hypotheses on the nearly \CY structure $(\oa, \Oa)$, the induced $G_2$-structure satisfies all the conditions in Theorem 3.9, and therefore can be deformed to have zero torsion. Finally, we pull back this torsion-free $G_2$-structure to obtain a genuine \CY structure on $M$. \\[0.5cm]

\renewcommand{\thethm}{3.10} 
\begin{thm}
\quad Let $\kappa > 0$ and $E_1, E_2, E_3, E_4 > 0$ be constants. Then there exist constants $\ep \in (0,1]$ and $K > 0$ such that whenever $0 < t \leq \ep$, the following is true. \\[-0.4cm]
\par Let $M$ be a compact 6-fold, and $(\oa, \Oa)$ a nearly \CY structure on $M$. Let $\oa'$, $g_M$ and $\ta_2'$ be as in $\S$3.1. Suppose
\begin{itemize}
	\item[\textnormal{(i)}] $\| \oa - \oa' \|_{L^2} \leq E_1 t^{3+\kappa},\ \,\| \oa - \oa' \|_{C^0} \leq E_1 t^{\kappa},\ \| \ta_2 - \ta_2' \|_{L^2} \leq E_1 t^{3+\kappa}$, \\[0.1cm]
and\ \,$\| \ta_2 - \ta_2' \|_{C^0} \leq E_1 t^{\kappa}$,
	\item[\textnormal{(ii)}] $\| \nab (\oa - \oa') \|_{L^{12}} \leq E_1 t^{-\frac{1}{2}+\kappa} $,\ \,$\| \nab \oa \|_{L^{12}} \leq E_1 t^{-\frac{1}{2}+\kappa}$, \\[0.1cm]
and\ \,$\| \nab \ta_1 \|_{L^{12}} \leq E_1 t^{-\frac{1}{2}+\kappa}$,  
	\item[\textnormal{(iii)}] $\| \nab (\oa - \oa') \|_{C^0} \leq E_2 t^{\kappa - 1}$,\ \,and\ \,$\| \nab^{2} (\oa - \oa') \|_{C^0} \leq E_2 t^{\kappa - 2}$,
	\item[\textnormal{(iv)}] the injectivity radius $\da(g_M)$ satisfies $\da(g_M) \geq E_3 t$, and
	\item[\textnormal{(v)}] the Riemann curvature $R(g_M)$ satisfies $\| R(g_M) \|_{C^0} \leq E_4 t^{-2}$. \\[-0.4cm]
 \end{itemize} 
Then there exists a \CY structure $(\tilde{J}, \tilde{\oa}, \tilde{\Oa})$ on $M$ such that $\| \tilde{\oa} - \oa \|_{C^0} \leq K t^{\kappa}$ and\, $\| \tilde{\Oa} - \Oa \|_{C^0} \leq K t^{\kappa}$. Moreover, if $H^1 (M, \BR) = 0$, then the cohomology classes satisfy $[\textnormal{Re}(\Oa)]$ = $[\textnormal{Re}(\tilde{\Oa})] \in H^3 (M, \BR)$\, and\, $[\oa]$ = $c\,[\tilde{\oa}] \in H^2 (M, \BR)$ for some $c > 0$. Here the connection $\nab$ and all norms are computed with respect to $g_M$. \\
\end{thm}
\pf \ \,Let $\vi$ be the 3-form on $\SM$ given by (3.11). Then Lemma 3.3 shows that $(\vi, g)$ is a $G_2$-structure on $\SM$, with $d\vi = 0$ as $d\oa = d\ta_1 = 0$. Define a 3-form $\psi = \vi - \ast_{g}\chi$ on $\SM$, where $\chi$ is the 4-form given by (3.11). Then $d^* \psi = d^* \vi - d^* (\ast_{g}\chi) = d^* \vi + \ast_{g} d\chi = d^* \vi$, since $d^* \ast_g = - \ast_g d$ on 3-forms and $d\chi = 0$. Now
$$ | \psi |_g = |\!\ast_g\!\psi |_g = |\!\ast_g\!\vi - \chi |_g \leq C\,|\!\ast_g\!\vi - \chi |_{g'} $$ 
where $C > 0$ is some constant relating norms w.r.t. the uniformly equivalent metrics $g$ and $g' = ds^2 + g_M$. From (3.15), one can show that $|\!\ast_g \vi - \chi |_{g'} \leq C_1 ( | \oa - \oa' |^{2}_{g_M} + | \oa - \oa' |_{g_M} + | \ta_2 - \ta_2' |_{g_M} ) $ for some $C_1 > 0$, and hence  
\renewcommand{\theequation}{3.16}
\begin{align}
| \psi |_g \leq C_2 \,(\,| \oa - \oa' |^{2}_{g_M} + | \oa - \oa' |_{g_M} + | \ta_2 - \ta_2' |_{g_M} )
\end{align}
for some $C_2 > 0$. Consequently, we have
$$\| \psi \|_{C^0} \leq C_2 \, (\,\| \oa - \oa' \|^{2}_{C^0} + \| \oa - \oa' \|_{C^0} + \| \ta_2 - \ta_2' \|_{C^0} ),$$
and
$$\| \psi \|_{L^2} \leq C_2 \, (\,\| \oa - \oa' \|_{C^0} \cdot \| \oa - \oa' \|_{L^2} + \| \oa - \oa' \|_{L^2} + \| \ta_2 - \ta_2' \|_{L^2} )$$
which then imply
\renewcommand{\theequation}{3.17}
\begin{align}
\| \psi \|_{C^0} \leq C_3 t^{\kappa}\quad \textnormal{and}\quad \| \psi \|_{L^2} \leq C_3 t^{3+\kappa}\ \,\textnormal{for some $C_3 > 0$},
\end{align}
where we have used condition (i) and $t \leq \ep \leq 1$. This verifies the first two inequalities of (i) in Theorem 3.9, and we now proceed to the last one. Denote by $\nab^g$ and $\nab^{g'}$ the connections computed using $g$ and $g'$ respectively. Since $d^* \psi = d^* \vi$, we get
\renewcommand{\theequation}{3.18}
\begin{align}
| d^* \psi |_g &= | d^* \vi |_g \leq | \nab^g \vi |_g \leq C\,| \nab^g \vi |_{g'}.   
\end{align}
Denote by $A$ the difference of the two torsion-free connections $\nab^g$ and $\nab^{g'}$. Thus $A$ transforms as  a tensor and it satisfies $A_{ij}^{k} = A_{ji}^{k}$. We need the following proposition to obtain the bound for $|d^* \psi|_g$. \\
\renewcommand{\theprop}{3.11}
\begin{prop}
In the situation above, we have 
\renewcommand{\theequation}{3.19}
\begin{align}
| \nab^{g'} \vi |_{g'} &\leq | \nab \oa |_{g_M} + | \nab \ta_1 |_{g_M}, \\ \tag{3.20}
| \nab^{g'} \vi' |_{g'} &\leq | \nab (\oa - \oa') |_{g_M} + | \nab \oa |_{g_M} + | \nab \ta_1 |_{g_M}, \\ \tag{3.21}
|A|_{g'} &\leq \frac{3}{2}\ |g^{-1}|_{g'} \cdot |\nab^{g'} g |_{g'}, \\ \tag{3.22}
and \ \  |\nab^{g'} g |_{g'} &\leq B_1\,| \nab^{g'}(\vi - \vi') |_{g'} + B_2\,| \vi - \vi' |_{g'} \cdot \big( | \nab^{g'} \vi |_{g'} + | \nab^{g'} \vi' |_{g'} \big) 
\end{align}
for some $B_1, B_2 > 0$ depending on a small upper bound for $| \vi - \vi' |_{g'}$. \\
\end{prop}
\pf \ \,For the first one, note that
\begin{align*} 
| \nab^{g'} \vi |_{g'} &= | \nab^{g'} (ds \we \oa + \ta_1) |_{g'} \notag \\
&\leq | \nab^{g'} ds |_{g'} \cdot | \oa |_{g'} + |ds|_{g'} \cdot | \nab^{g'} \oa |_{g'} + | \nab^{g'} \ta_1 |_{g'}. \notag 
\end{align*}
Since $ds$ is a constant 1-form with length 1 w.r.t. the metric $g'$, equation (3.19) follows. The second inequality follows easily from the first one. For (3.21), we have from the definition of the tensor $A$,
$$ A_{ij}^{k} = \Ga_{ij}^{k} - \Ga_{ij}^{'k} $$
where $\Ga_{ij}^{k}$ and $\Ga_{ij}^{'k}$ are the Christoffel symbols of the Levi-Civita connections $\nab^g$ and $\nab^{g'}$ respectively. Consider now the term $\nab^{g'}g$, and expressing it in index notation,
$$ \nab^{g'}_{a} g_{bc} = \pa_a g_{bc} - \Ga_{ab}^{'d}g_{dc} - \Ga_{ac}^{'d}g_{bd}. $$
Then 
\begin{align*}
\Ga_{ij}^{k} &= \frac{1}{2}\,g^{kl}\big( \pa_i g_{jl} + \pa_j g_{il} - \pa_l g_{ij} \big) \\
&= \frac{1}{2} g^{kl} \big[ \big( \nab^{g'}_i g_{jl} + \Ga_{ij}^{'m}g_{ml} + \Ga_{il}^{'m}g_{jm} \big) + \big( \nab^{g'}_j g_{il} + \Ga_{ji}^{'m}g_{ml} + \Ga_{jl}^{'m}g_{im} \big) \\
&\quad - \big( \nab^{g'}_l g_{ij} + \Ga_{li}^{'m}g_{mj} + \Ga_{lj}^{'m}g_{im} \big) \big] \\
&= \frac{1}{2}\,g^{kl} \big( \nab^{g'}_i g_{jl} + \nab^{g'}_j g_{il} - \nab^{g'}_l g_{ij} + 2\,\Ga_{ij}^{'m}g_{ml} \big) \\
\end{align*} 
by the fact that $\nab^{g'}$ is torsion-free. Hence,
\begin{align*}
A_{ij}^{k} &= \Ga_{ij}^{k} - \Ga_{ij}^{'k} \\
&= \frac{1}{2}\,g^{kl} \big( \nab^{g'}_i g_{jl} + \nab^{g'}_j g_{il} - \nab^{g'}_l g_{ij} \big),
\end{align*}
which gives rise to (3.21). \\ 

To verify the last inequality, first let $F$ be the smooth function that maps each positive 3-form to its associated metric. Then $F(\vi) = g$ and $F(\vi') = g'$. As $F(\vi)$ depends pointwise on $\vi$, we can write $F(\vi)(x) = F(x, \vi(x))$ for all $x \in \SM$ and $\vi (x)$ in the vector space $\bigwedge^{3} T^{*}_{x} (\SM)$. We may then take partial derivative in the $\vi(x)$ direction without using a connection, and write $\pa$ for the partial derivative in this direction. Now,
\begin{align*}
|\nab^{g'} g |_{g'} &= |\nab^{g'} (g - g') |_{g'} \\
&= |\nab^{g'} (F(\vi) - F(\vi')) |_{g'} \\
&= \Big| \int_{0}^{1} \frac{d}{dr} \nab^{g'} \big( F(x, \vi'(x) + r(\vi(x) - \vi'(x))) \big)\,dr  \Big|_{g'} \\
&= \Big| \int_{0}^{1} \nab^{g'} \big( \frac{d}{dr} F(x, \vi'(x) + r(\vi(x) - \vi'(x))) \big)\,dr \Big|_{g'} \\
&= \Big| \int_{0}^{1} \nab^{g'} \big( \pa F(x, \vi'(x) + r(\vi(x) - \vi'(x))) \cdot (\vi(x) - \vi'(x)) \big)\,dr \Big|_{g'} \\
& = \Big| \int_{0}^{1} \big[ (\nab^{g'} \pa F) (x, \vi'(x) + r(\vi(x) - \vi'(x)))  \\
&\qquad + \pa^2 F(x, \vi'(x) + r(\vi(x) - \vi'(x))) \cdot \nab^{g'}(\vi'(x) + r(\vi(x) - \vi'(x))) \big]  \\
&\qquad \cdot (\vi(x) - \vi'(x)) + \pa F(x, \vi'(x) + r(\vi(x) - \vi'(x))) \\
&\qquad \cdot \nab^{g'}(\vi(x) - \vi'(x))\,dr \Big|_{g'}. 
\end{align*} 
It can be shown, by using the fact that continuous functions are bounded over compact spaces, that for any $\phi$ which is close enough to $\vi'$, we have
$$|\pa F(x, \phi(x))|_{g'} \leq B_1  \quad \textnormal{and}\quad |\pa^2 F(x, \phi(x))|_{g'} \leq 2B_2 $$ 
for some constants $B_1, B_2 > 0$, and as this is a calculation at a point, $B_1, B_2$ are constants depend only on a small upper bound for $|\vi - \vi'|_{g'}$. Moreover, if we choose geodesic normal coordinates at $x$, then the Christoffel symbols $\Ga_{ij}^{'k}$ of $\nab^{g'}$ vanish at $x$, so $\nab^{g'}$ reduces to the usual partial differentiation at $x$ and it follows that 
$$\nab^{g'} \pa F (x, \vi'(x) + r(\vi(x) - \vi'(x))) = 0 $$
since $F$, and hence $\pa F$ is invariant under translation along the directions of coordinate vectors. Consequently,   
\begin{align*}
|\nab^{g'} g |_{g'} &\leq B_1 |\nab^{g'}(\vi(x) - \vi'(x))|_{g'} \\ 
&\quad + 2\, B_2 \,|\vi(x) - \vi'(x)|_{g'} \cdot \int_{0}^{1} |\nab^{g'} (\vi'(x) + r(\vi(x) - \vi'(x))) |_{g'}\, dr \\
&\leq B_1 |\nab^{g'}(\vi(x) - \vi'(x))|_{g'} + B_2 \,|\vi(x) - \vi'(x)|_{g'}\!\cdot\!\big( | \nab^{g'} \vi |_{g'} + | \nab^{g'} \vi' |_{g'} \big) 
\end{align*}
and this finishes the proof of Proposition 3.11.  \hfill $\Box$ \\

\par Applying the above estimates to (3.18) shows that  
\begin{align*}
| d^* \psi |_{g} &\leq C |\nab^g \vi |_{g'} \\
&\leq C |\nab^{g'} \vi |_{g'} + C |A|_{g'} \cdot |\vi|_{g'} \\
&\leq C | \nab \oa |_{g_M} + C | \nab \ta_1 |_{g_M} + \frac{3}{2} \, C \ |g^{-1}|_{g'} \cdot |\nab^{g'} g |_{g'} \cdot |\vi|_{g'} \\
&\qquad \textnormal{by (3.19) and (3.21)} \\[0.1cm]
&\leq C | \nab \oa |_{g_M} + C | \nab \ta_1 |_{g_M} + C_4 \big( B_1\,| \nab^{g'}(\vi - \vi') |_{g'} \\
&\quad + B_2\,| \vi - \vi' |_{g'} \cdot ( | \nab^{g'} \vi |_{g'} + | \nab^{g'} \vi' |_{g'} )  \big) \\
&\qquad \textnormal{by (3.22) and the fact that $g^{-1}$ and $\vi$ are bounded by some}\\
&\qquad \textnormal{constants w.r.t. $g'$} \\[0.1cm]
&\leq C | \nab \oa |_{g_M} + C | \nab \ta_1 |_{g_M} + C_5 |\nab (\oa - \oa') |_{g_M} + C_6  |\oa - \oa' |_{g_M} \cdot \\
&\qquad (| \nab (\oa - \oa') |_{g_M} + 2 | \nab \oa |_{g_M} + 2 | \nab \ta_1 |_{g_M}) \quad \textnormal{by (3.20),}
\end{align*}
where $C_4 , C_5 , C_6 > 0$ are some constants. From the second inequality of (i), we have $\|\oa - \oa' \|_{C^0} \leq \lda t^{\kappa} \leq \lda $ as $t \leq \kappa \leq 1$. Combining with condition (ii) shows that 
\begin{align*}
\| d^* \psi \|_{L^{12}} &\leq C \| \nab \oa \|_{L^{12}} + C \| \nab \ta_1 \|_{L^{12}} + C_5 \|\nab (\oa - \oa') \|_{L^{12}} \\ 
&\ \, + C_6 \|\oa - \oa' \|_{C^0} \cdot ( \|\nab (\oa - \oa') \|_{L^{12}} + 2 \| \nab \oa \|_{L^{12}} + 2 \| \nab \ta_1 \|_{L^{12}}) \\
&\leq (C\lda + C\lda + C_5 \lda + C_6 \lda \,(\lda + 2 \lda + 2 \lda)) t^{-\frac{1}{2}+\kappa}. 
\end{align*}
Thus, together with (3.17), we have verified condition (i) of Theorem 3.9. \\

\par Given (iv) and (v), the injectivity radius and the Riemann curvature of $g' = ds^2 + g_M$ satisfy $\da(g') \geq$ min($\mu t, \pi$) = $\mu t$ for small $t$, and $\| R(g') \|_{C^0} \leq \nu t^{-2}$. Basically, the role of the estimates on injectivity radius and the $C^0$-norm of the Riemann curvature in the proof of Theorem 3.9 is to show that there exist coordinate systems on small balls such that the metric in these coordinate systems is $C^{1,\al}$-close to the Euclidean metric on $\BR^7$ for $\al \in (0,1)$. Therefore, all we have to show is that the metric $g$ is $C^{1,\al}$-close to the Euclidean metric. Now we know that $g'$ is $C^{1,\al}$-close to the Euclidean metric on $\BR^7$ by Joyce \cite[Prop. 11.7.2]{Joyce1}. But the second inequality of condition (i) and condition (iii) of the hypotheses ensure that $g$ and $g'$ are $C^2$-close, and this implies that $g$ is also $C^{1,\al}$-close to the Euclidean metric on $\BR^7$, which is what we need. \\

\par Therefore Theorem 3.9 gives a torsion-free $G_2$-structure $(\tilde{\vi}, \tilde{g})$ on $\SM$. It remains to construct a \CY structure on $M$ from $(\tilde{\vi}, \tilde{g})$. Denote $\frac{\pa}{\pa s}$ by the Killing vector w.r.t. $g$ such that $\iota \left(\frac{\pa}{\pa s}\right)ds = 1$. Then $\frac{\pa}{\pa s}$ is also a Killing vector w.r.t. $\tilde{g}$ since $(\tilde{\vi}, \tilde{g})$ is $S^1$-invariant. It follows from a general fact about Killing vectors on a torsion-free, compact $G_2$-manifold that $\nab^{\tilde{g}} \frac{\pa}{\pa s} = 0$, and hence $\left|\frac{\pa}{\pa s}\right|_{\tilde{g}}$ equals to some constant $c$. Define a 1-form $d\tilde{s}$ on $\SM$ by $(d\tilde{s})_a = \frac{1}{c}\,\tilde{g}_{ab}\big(\frac{\pa}{\pa s}\big)^b $. Then $d\tilde{s}$ is closed, of unit length w.r.t. $\tilde{g}$, and $\iota\left(\frac{\pa}{\pa s}\right)d\tilde{s} = c$, and we may write $d\tilde{s} = c\,ds + \al'$ for some closed 1-form $\al'$ on $\SM$ with $\iota\left(\frac{\pa}{\pa s}\right) \al' = 0$. The 1-form $\al'$ is thus the pullback of some closed 1-form $\al$ on $M$ via the projection map $\pi : \SM \longrightarrow M$, i.e. $\al' = \pi^* (\al)$. Since by assumption $H^1 (M, \BR) = 0$, $\al$, and hence $\al'$, is exact. Therefore we have $[d\tilde{s}] = c\,[ds]$. \\

\par Using the fact that $\frac{\pa}{\pa s}$ is a Killing vector and $\tilde{\vi}$ is a closed 3-form, we have
$$ d \Big( \iota\Big(\frac{\pa}{\pa s}\Big) \tilde{\vi} \Big) = \mathcal{L}_{\frac{\pa}{\pa s}} \tilde{\vi} = 0, $$
so $\iota\,\big(\frac{\pa}{\pa s}\big)\,\tilde{\vi}$, and similarly $\iota\,\big(\frac{\pa}{\pa s}\big)\,(\ast_{\tilde{g}}\tilde{\vi})$ are both $S^1$-invariant closed forms on $\SM$. Then we can define a closed 2-form $\tilde{\oa}$ and closed 3-forms $\tilde{\ta_1}$ and $\tilde{\ta_2}$ on $M$ by 
\begin{align*}
\tilde{\oa}_x = \frac{1}{c}\,\Big( \iota\,\Big(\frac{\pa}{\pa s}\Big)\,\tilde{\vi} \Big)_{(s,x)} \Big|_{T_x M}, \quad (\tilde{\ta_1})_x = \tilde{\vi}_{(s,x)} - (d\tilde{s} \we \tilde{\oa})_{(s,x)} |_{T_x M}, \\[0.3cm]
\textnormal{and }\quad (\tilde{\ta_2})_x = -\frac{1}{c}\,\Big( \iota\,\Big(\frac{\pa}{\pa s}\Big)\,\big(\!\ast_{\tilde{g}}\tilde{\vi}\big)\Big)_{(s,x)} \Big|_{T_x M}
\end{align*}
for each $x \in M$ and any $s \in S^1$. Identify $T_{(s,x)} (\SM)$ with $\BR^7 \cong \BR \oplus \BC^3$ such that $T_x M$ is identified with $\BC^3$, $(\tilde{\vi}, \tilde{g})$ with the flat $G_2$-structure $(\vi_0, g_0)$, $\frac{1}{c}(\frac{\pa}{\pa s})$ with $\frac{\pa}{\pa x_1}$, and $d\tilde{s}$ with $dx_1$ where $x_1$ is the coordinate on $\BR$. Then calculation shows that at each $x \in M$, $(\tilde{J}, \tilde{\oa}, \tilde{\Oa}) $ can be identified with the standard \CY structure $(J_0, \oa_0, \Oa_0)$ on $\BC^3$, where $\tilde{\Oa} = \tilde{\ta_1} + i \tilde{\ta_2}$, and $\tilde{J}$ is the associated complex structure. It follows that $(\tilde{J}, \tilde{\oa}, \tilde{\Oa} )$ gives a \CY structure on $M$  with $\tilde{\vi} = d\tilde{s} \we \tilde{\oa} + \tilde{\ta_1}$\, and\, $\ast_{\tilde{g}} \tilde{\vi} = \frac{1}{2} \tilde{\oa} \we \tilde{\oa} - d\tilde{s} \we \tilde{\ta_2}$ on $\SM$.   \\



\par It is not hard to show $\| \tilde{\oa} - \oa \|_{C^0} \leq K t^{\kappa}$ and, by making $K$ larger if necessary, $\| \tilde{\Oa} - \Oa \|_{C^0} \leq K t^{\kappa}$ provided that $\| \tilde{\vi} - \vi \|_{C^0} \leq K t^{\kappa}$, which is a consequence from Theorem 3.9. Moreover, as $[\tilde{\vi}] = [\vi]$ and $[d\tilde{s}] = c\,[ds]$, it follows that $[\oa] = c\,[\tilde{\oa}]$ and $[\tilde{\ta}_1] = [\ta_1]$. This completes the proof of Theorem 3.10. \hfill $\Box$ \\[1cm]  
\bf Remarks
 \normalfont
\begin{enumerate}
	\item In general we can't guarantee [Im($\tilde{\Oa}$)] = [Im($\Oa$)] since, roughly speaking, [Im($\tilde{\Oa}$)] is locally determined by [Re($\tilde{\Oa}$)], whereas [Im($\Oa$)] is free to change slightly, as long as the inequality (3.2) is satisfied. Hence [Im($\Oa$)] is independent of [Re($\Oa$)] and it follows that [Im($\tilde{\Oa}$)] can't possibly be determined by [Im($\Oa$)]. \\
	\item If $H^1 (M, \BR) \ne 0$, then $\al'$ may not be exact, and we have to modify the cohomological formula for [Re($\tilde{\Oa}$)] to 
	$$ [\textnormal{Re}(\tilde{\Oa})] = [\textnormal{Re}(\Oa)] - [\al] \cup [\tilde{\oa}].$$ 
	\item There is an alternative way of obtaining the \CY structure on $M$ from the holonomy point of view. Since $(\tilde{\vi}, \tilde{g})$ is torsion-free, Hol($\tilde{g}$) $\subseteq G_2$. Moreover, Hol($\tilde{g}$) fixes the vector $\frac{\pa}{\pa s}$ as $\nab^{\tilde{g}} \frac{\pa}{\pa s} = 0$. It turns out that Hol($\tilde{g}$) actually lies in SU(3) and hence the torsion-free $G_2$-structure $(\tilde{\vi}, \tilde{g})$ must be coming from a \CY structure on $M$.  
\end{enumerate}

\begin{center}
\large \bf{4. \ Calabi--Yau cones, Calabi--Yau manifolds with a conical singularity and Asymptotically Conical Calabi--Yau manifolds} \\[0.5cm]
\end{center}





\normalsize

\par In this section we define \CY cones, \CY manifolds with a conical singularity and Asymptotically Conical \CY manifolds. We will give some examples and provide results analogous to the usual Darboux Theorem on symplectic manifolds for the \CY manifolds with a conical singularity and Asymptotically Conical \CY manifolds. The conical singularities in \CY 3-folds will be desingularized in section 5 by using the existence result obtained in section 3. \\

\begin{center}
\bf{4.1. \ Preliminaries on Calabi--Yau cones}\\[0.3cm]
\end{center}

\par We will give our definition of \CY cones and provide several examples in this section. Let us first consider the $\BC^m$ case. Write $\BC^m$ as $S^{2m-1} \times (0,\infty) \cup \{0\}$, a cone over the $(2m-1)$-dimensional sphere. Let $r$ be a coordinate on $(0, \infty)$. Then the standard metric $g_0$, \ka form $\oa_0$ and holomorphic volume form $\Oa_0$ on $\BC^m$ can be written as 
\begin{align*}
g_0 = r^2 g_0 &|_{S^{2m-1}} + dr^2, \quad \oa_0 = r^2 \oa_0 |_{S^{2m-1}} + rdr \we \al \\[0.1cm]
&\textnormal{and}\quad \Oa_0 = r^m \Oa_0 |_{S^{2m-1}} + r^{m-1}dr \we \ba,
\end{align*} 
where $\al$ is a real 1-form and $\ba$ a complex $(m-1)$-form on $S^{2m-1}$. Hence they scale as
\begin{align*}
g_0 |_{S^{2m-1} \times \{r\}}& = r^2 g_0 |_{S^{2m-1}}, \quad \oa_0 |_{S^{2m-1} \times \{r\}} = r^2 \oa_0 |_{S^{2m-1}} \\[0.15cm]
&\textnormal{and}\quad \Oa_0 |_{S^{2m-1} \times \{r\}} = r^m \Oa_0 |_{S^{2m-1}}.
\end{align*} 
Motivated by this standard case, we give our definition of a \emph{\CY cone}: \\

\renewcommand{\thedefinition}{4.1}
\begin{definition}
\quad \textnormal{Let $\Gamma$ be a compact ($2m-1$)-dimensional smooth manifold, and let $V = \{0\} \cup V'$ where $V' = \Gamma \times (0, \infty)$. Write points on $V'$ as $(\gamma, r)$. $V$ is called a \emph{\CY cone} if $V'$ is a \CY $m$-fold with a \CY structure $(J_V , \oa_V , \Oa_V)$ and its associated \CY metric $g_V$ satisfying 
\renewcommand{\theequation}{4.1}
\begin{align}
g_V = &r^2 g_V |_{\Ga \times \{1\}} + dr^2, \quad \oa_V = r^2 \oa_V |_{\Ga \times \{1\}} + rdr \we \al \notag \\[-0.2cm]
\\[-0.2cm]
&\textnormal{and}\quad \Oa_V = r^m \Oa_V |_{\Ga \times \{1\}} + r^{m-1}dr \we \ba. \notag
\end{align}
Here we identify $\Ga$ with $\Ga \times \{1\}$, and $\al$ is a real 1-form and $\ba$ a complex $(m-1)$-form on $\Ga$.
} \\
\end{definition} 

\par Let $X$ be the radial vector field on $V$ such that $X_{(\ga,r)}= r\frac{\pa}{\pa r}$ for any $(\ga,r) \in \Ga \times (0, \infty)$. Then $r^2 \al = \iota(X)\oa_V$ and $r^m \ba = \iota(X)\Oa_V$. Moreover,
\begin{align*}
\mathcal{L}_X \oa_V &= d(\iota(X)\oa_V) \quad \textnormal{as $d\oa_V = 0$}\\
&= d(r^2 \al)\\
&= r^2 d\al + 2rdr \we \al. 
\end{align*}
It can be shown that $d\al = 2\oa_V |_{\Ga}$ by using $d\oa_V = 0$ and the formula for $\oa_V$ in (4.1). Therefore we have $\mathcal{L}_X \oa_V = 2\oa_V$. In a similar way, we can show $\mathcal{L}_X \Oa_V = m \Oa_V$ and $\mathcal{L}_X g_V = 2g_V$. It follows that $\mathcal{L}_X J_V = 0$, and hence $X$ is a holomorphic vector field. The flow of $X$ thus expands the \CY metric $g_V$, the \ka form $\oa_V$ and the holomorphic volume form $\Oa_V$ exponentially and $X$ is then a Liouville vector field. In particular, the 1-form $\al$ defines a contact form on $\Ga$, which makes $\Ga$ a contact $(2m-1)$-fold. \\
 
\par The tangent space $T_{(\ga,r)} V$ decomposes as $T_{(\ga,r)} V = T_{\ga} \Ga \, \oplus \! <\!X_{(\ga,r)}\!>_{\BR}$ for any $(\ga, r) \in \Ga \times (0, \infty)$. Note that $Z := J_V X$ is a vector field on $\Ga$, and it is complete as $\Ga$ is compact. Now $\iota(Z)\oa_V$ is a 1-form such that $\iota(X)(\iota(Z)\oa_V) = g_V (X, X) = r^2$ and $\iota(Z)\oa_V \big|_{\Ga \times \{r\}} = 0$, hence we can write $\iota(Z)\oa_V = rdr$. It follows that
$$\mathcal{L}_{Z} \oa_V = d(\iota(Z)\oa_V) = d(rdr) = 0.   $$ \\[-0.2cm]
For the holomorphic volume form, we use the fact that if $\Oa$ is a holomorphic $(m,0)$-form and $v$ a holomorphic vector field, then $\mathcal{L}_{Jv} \Oa = i \mathcal{L}_{v} \Oa $ where $J$ is the complex structure. Now $Z = J_V X$ is a holomorphic vector field, this gives $\mathcal{L}_{Z} \Oa_V = im \Oa_V$. \\

\par Now we define a \emph{complex dilation} on the \CY cone $V$. The flow of $Z$ generates the diffeomorphism $\textnormal{exp} (\ta Z)$ on $\Ga$ for each $\ta \in \BR$. Thus for each $\ta \in \BR$ and $t>0$, we can define a \emph{complex dilation} $\lda$ on $V$ which is given by $\lda(0) = 0$ and $\lda(\ga, r) = (\textnormal{exp} (\ta Z)(\ga),\, tr)$. \\



\renewcommand{\thelem}{4.2} 
\begin{lem}
\quad Let $\lda : V \longrightarrow V$ be the complex dilation defined above. Then $\lda^* (g_V) = t^2 g_V$, $\lda^* (\oa_V) = t^2 \oa_V$ and $\lda^* (\Oa_V) = t^m e^{im\ta}\,\Oa_V$. \\
\end{lem}
\pf \ \, It follows from $\mathcal{L}_{Z} \oa_V = 0$ that $\textnormal{exp}(\ta Z)^* (\oa_V) = \oa_V$ and hence $\lda^* (\oa_V) = t^2 \oa_V$ by the scaling of $t$. The formula for the metric $g_V$ follows similarly. For the holomorphic $(m,0)$-form $\Oa_V$, observe that
\begin{align*}
\frac{d}{d\ta}\,\textnormal{exp}(\ta Z)^* (\Oa_V)\big|_{\ta = 0} = \mathcal{L}_{Z} \Oa_V = im\,\Oa_V,
\end{align*}
and this means $\textnormal{exp}(\ta Z)^* (\Oa_V) = e^{im\ta}\,\Oa_V$. Thus together with the scaling of $t$, we have $\lda^* (\Oa_V) = t^m e^{im\ta}\,\Oa_V$. \hfill $\Box$ \\

\par In the situation of our standard example $\BC^m$, the complex dilation is given by complex multiplication $\lda : \BC^m \longrightarrow \BC^m$ sending $z$ to $\lda z$, where $\lda = te^{i\ta} \in \BC$. It is easy to see the above properties for the standard structures $g_0$, $\oa_0$ and $\Oa_0$. \\[0.6cm]
\textbf{Example 4.3.}\quad A trivial example is given by $\BC^m$, a cone on $S^{2m-1}$. Some nontrivial examples can be constructed as follows: Let $G$ be a finite subgroup of SU($m$) acting freely on $\BC^m \setminus \{0\}$, then the quotient singularity $\BC^m /G$ is a \CY cone. An example of this type is given by the $\mathbb{Z}_m$-action on $\BC^m$: 
\begin{align*}
\zeta^k \cdot (z_1 , \dots , z_m) \ = \ (\zeta^k \,z_1 , \dots , \zeta^k \,z_m)
\end{align*}
where $\zeta = e^{2 \pi i / m}$ and $1 \leq k \leq m$. Note that $\zeta^m = 1$, so $\mathbb{Z}_m$ is a subgroup of SU($m$) and acts freely on $\BC^m \setminus\!\{0\}$. Then $\BC^m  / \mathbb{Z}_m$ is a \CY cone. \\[0.6cm]
\textbf{Example 4.4.}\quad Consider the cone $V$ defined by the quadric $\sum^{m+1}_{j=1} z_{j}^2 = 0$ on $\BC^{m+1}$. The singularity at the origin is known as an \emph{ordinary double point}, or a \emph{node}. 
It can be shown that $\Ga$ is an $S^{m-1}$-bundle over $S^{m}$. Stenzel \cite{Ste} constructed a \CY metric on $V$, thus making it a \CY cone. We are particularly interested in the case $m =3$. Then $\Ga$ has the topology of $S^2 \times S^3$, and hence $V$ is topologically a cone on $S^2 \times S^3$ since any $S^2$-bundle over $S^3$ is trivial. One can also describe $V$ as follows. Consider the blow-up $\widetilde{\BC}^4$ of $\BC^4$ at origin. It introduces an exceptional divisor $\BC \BP^3$, and the blow-up $\widetilde{V}$ of the cone $V$ inside $\widetilde{\BC}^4$ meets this $\BC \BP^3$ at $S \cong \BC \BP^1 \times \BC \BP^1$. The exceptional divisor $\BC \BP^3$ corresponds to the zero section of the line bundle $L$ given by $\widetilde{\BC}^4 \longrightarrow \BC \BP^3$, and so its normal bundle is isomorphic to $L$. Hence the normal bundle $\mathcal{O}(-1,-1)$ over $\BC \BP^1 \times \BC \BP^1$ is isomorphic to the line bundle $\widetilde{V} \longrightarrow S$. This gives us the following isomorphisms: 
$$ V \setminus \{0\} \ \cong \ \widetilde{V} \setminus S \ \cong \ \mathcal{O}(-1,-1) \setminus (\BC \BP^1 \times \BC \BP^1). $$ \\[0.8cm]
\textbf{Example 4.5.}\quad Suppose $S$ is K\"{a}hler-Einstein with positive scalar curvature, Calabi \cite[p.284-5]{Calabi} constructed a 1-parameter family of \CY metrics $g_t$ for $t \geq 0$ on the canonical line bundle $K_S$. When $t > 0$, $g_t$ is a nonsingular complete metric on $K_S$ and when $t = 0$, $g_0$ degenerates on $S$ and thus gives a cone metric on $K_S \setminus S$, which then makes $K_S \setminus S$ a \CY cone with $S$ ``collapsed" to the vertex of the cone. \\
\begin{itemize}
		\item[(i)] One of the standard examples of K\"{a}hler-Einstein manifolds with positive scalar curvature is given by the complex surface $S \cong \BC \BP^1 \times \BC \BP^1$. Calabi's construction thus applies to it and yields a \CY metric on $K_{S} = \mathcal{O}(-2,-2) \longrightarrow S$. Note that $\mathcal{O}(-1,-1)$ is a double cover of $\mathcal{O}(-2,-2)$ away from the zero section $S$, so we have the following relation between the cone $K_S \setminus S$ and the cone $V$ described in Example 4.4:
$$K_S \setminus S \ \cong \ (V \setminus \{0\})/\mathbb{Z}_2 .$$	\\[-0.5cm]
		\item[(ii)] Boyer, Galicki and Koll$\acute{\textnormal{a}}$r \cite{BGK} constructed K\"{a}hler-Einstein metrics on some compact orbifolds, particularly on orbifolds of the form given by the quotient $(L \setminus \{0\}) / \BC^*$, where $L = \{ (z_1, \dots, z_m) \in \BC^m : \sum_{j=1}^{m} z_{j}^{a_j} = 0 \}$ for some positive integers $a_j$ satisfying certain conditions. The set $L$ is a hypersurface in $\BC^m$, and $\BC^*$ acts naturally on $\BC^m$ by $\lda : (z_1, \dots , z_m) \mapsto (\lda^{a_1}z_1, \dots, \lda^{a_m}z_m) $. It follows from Calabi's construction, in the category of orbifolds, that $L \setminus \{0\}$ is a \CY cone. \\[1.5cm]
\end{itemize}

\begin{center}
\bf {4.2. \ \CY \bf{\emph{m}}-folds with a conical singularity} \\[0.3cm]
\end{center}

\normalsize

\par We define \CY $m$-folds with a conical singularity in this section. Moreover, we show that there exist coordinate systems that can trivialize the symplectic forms of \CY $m$-folds with a conical singularity. \\

\renewcommand{\thedefinition}{4.6}
\begin{definition}
\quad \textnormal{Let $(M_0, J_0, \oa_0, \Oa_0)$ be a singular \CY $m$-fold with a singularity at $x \in M_0$, and no other singularities. We say that $M_0$ is a \emph{\CY $m$-fold with a conical singularity at $x$} with rate $\nu > 0$ modelled on a \CY cone $(V, J_V , \oa_V , \Oa_V)$ if there exist a small $\ep > 0$, a small open neighbourhood $S$ of $x$ in $M_0$, and a diffeomorphism $\Phi : \Ga \times (0, \ep) \longrightarrow S \setminus \{x\}$ such that 
\renewcommand{\theequation}{4.2}
\begin{align}
|\nab^k ( \Phi^* ( \oa_0 ) - \oa_V ) |_{g_V} &= O(r^{\nu - k}), \ \textnormal{and} \\ \tag{4.3}
|\nab^k ( \Phi^* ( \Oa_0 ) - \Oa_V ) |_{g_V} &= O(r^{\nu - k}) \ \ \textnormal{as } r \rightarrow 0 \ \textnormal{and for all } k \geq 0. 
\end{align}
Here $\nab$ and $|\cdot|$ are computed using the cone metric $g_V$.} \\
\end{definition}

\par Note that the asymptotic conditions on $g_0$ and $J_0$ follow from those on $\oa_0$ and $\Oa_0$, namely,
\renewcommand{\theequation}{4.4}
\begin{align}
|\nab^k ( \Phi^* ( g_0 ) - g_V ) |_{g_V} &= O(r^{\nu - k}), \ \textnormal{and} \\ \tag{4.5}
|\nab^k ( \Phi^* ( J_0 ) - J_V ) |_{g_V} &= O(r^{\nu - k}) \ \ \textnormal{as } r \rightarrow 0 \ \textnormal{and for all } k \geq 0, 
\end{align}
and so it is enough to just assume asymptotic conditions on $\oa_0$ and $\Oa_0$. \\

\par We will usually assume that $M_0$ is compact. The point of the definition is that $M_0$ has one end modelled on $\Ga \times (0, \ep)$ and as $r \rightarrow 0$, all the structures $g_0, J_0, \oa_0$ and $\Oa_0$ on $M_0$ converge to the cone structures $g_V, J_V, \oa_V$ and $\Oa_V$ with rate $\nu$ and with all their derivatives. \\

\par The 2-forms $\Phi^* ( \oa_0 )$ and $\oa_V$ are closed on $\Ga \times (0,\ep)$ and so $\Phi^* ( \oa_0 ) - \oa_V$ represents a cohomology class in $H^2 (\Ga \times (0,\ep), \BR) \cong H^2 (\Ga, \BR)$. Similarly, $\Phi^* ( \Oa_0 ) - \Oa_V$ represents a cohomology class in $H^m (\Ga \times (0,\ep), \BC) \cong H^m (\Ga, \BC)$. It turns out that in the conical singularity case, these two classes $[\Phi^* ( \oa_0 ) - \oa_V]$ and $[\Phi^* ( \Oa_0 ) - \Oa_V]$ are automatically zero. \\

\renewcommand{\thelem}{4.7}
\begin{lem}
\quad Let $(M_0, J_0, \oa_0, \Oa_0)$ be a compact \CY $m$-fold with a conical singularity at $x$ with rate $\nu >0$ modelled on a \CY cone $(V, J_V , \oa_V , \Oa_V)$. Then $[\Phi^* ( \oa_0 ) - \oa_V ] = 0$ in $H^2 (\Ga, \BR)$ and $[\Phi^* ( \Oa_0 ) - \Oa_V ] = 0$ in $H^m (\Ga, \BC)$. \\
\end{lem}
\pf \ \,Suppose $\Si$ is a 2-cycle in $\Ga$. Then
\begin{align*} 
\int_{\Si \times \{r\}} \big( \Phi^* ( \oa_0 ) - \oa_V \big) &= \textnormal{vol}(\Si) \cdot O(r^{\nu}) \quad \textnormal{by (4.2)} \\[-0.2cm]
&= O(r^{\nu+2})
\end{align*}
Hence $\nu > 0$ implies the above integral approaches 0 as $r \rightarrow 0$. Then 
$$ [\Phi^* ( \oa_0 ) - \oa_V ] \cdot [\Si] = 0 $$
for any 2-cycle $\Si$, and hence $[\Phi^* ( \oa_0 ) - \oa_V] = 0 \in H^2 (\Ga \times (0, \ep), \BR) \cong H^2 (\Ga, \BR)$. The case for $[\Phi^* ( \Oa_0 ) - \Oa_V ]$ follows similarly by considering 3-cycles in $\Ga$. \hfill $\Box$ \\

\par One may ask whether the symplectic form $\oa_0$ on $S$ near $x$ in $M_0$ can actually be symplectomorphic to the cone form $\oa_V$ near the origin, rather than just having an asymptotic relation in (4.2). Our next result shows that this is indeed the case for \CY $m$-folds with conical singularities, and can be regarded as an analogue of the usual Darboux theorem on symplectic manifolds.   \\

\renewcommand{\thethm}{4.8}
\begin{thm}
\quad Let $(M_0, J_0, \oa_0, \Oa_0)$ be a compact \CY $m$-fold with a conical singularity at $x$ with rate $\nu >0$ modelled on a \CY cone $(V, J_V , \oa_V , \Oa_V)$. Then there exist an $\ep' > 0$, an open neighbourhood $S'$ of $x$ and a diffeomorphism $\Phi' : \Ga \times (0, \ep') \longrightarrow S' \setminus \{x\}$ such that $(\Phi')^* (\oa_0) = \oa_ V$ and $|\nab^k ( (\Phi')^* ( \Oa_0 ) - \Oa_V ) |_{g_V} = O(r^{\nu - k})$ for all $k \geq 0$.  \\
\end{thm}
\pf \ \,By the definition of a \CY $m$-fold with a conical singularity, there is an $\ep > 0$ and a diffeomorphism $\Phi : \Ga \times (0, \ep) \longrightarrow S \setminus \{x\}$ such that $|\nab^k ( \Phi^* ( \oa_0 ) - \oa_V ) |_{g_V} = O(r^{\nu - k})$ as $r \rightarrow 0$ and for all $k \geq 0$. Write $\oa = \Phi^* (\oa_0)$. Thus we have two symplectic forms, namely the cone form $\oa_V$ and the asymptotic cone form $\oa$ on $\Ga \times (0, \ep)$. Following Moser's proof of the usual Darboux Theorem \cite[p.93]{McDuff}, we construct a 1-parameter family of closed 2-forms \\[-0.5cm]
\renewcommand{\theequation}{4.6}
\begin{align}
\oa^t = \oa_V + t(\oa - \oa_V) \quad \textnormal{for } t \in [0,1] \\[-0.5cm] \notag
\end{align}
on $\Ga \times (0, \ep)$. Make $\ep$ smaller if necessary so that $\oa^t$ is symplectic for all $t \in [0,1]$. Then 
\renewcommand{\theequation}{4.7}
\begin{align}
\frac{d}{dt} \oa^t = \oa - \oa_V. \\[-0.3cm] \notag
\end{align}
By Lemma 4.7, $[\oa - \oa_V] = 0$ in $H^2 (\Ga, \BR)$ and hence $\oa - \oa_V$ is exact. 
Suppose $\eta = \oa - \oa_V$, and write $\eta$ as $\eta_0(\ga,r) + \eta_1 (\ga,r) \we dr$, where $\eta_0(\ga,r) \in \Lambda^2 T^{*}_{\ga} \Ga$ and $\eta_1(\ga,r) \in \Lambda^1 T^{*}_{\ga} \Ga$. Then $\eta$ is an exact 2-form such that
$$ \frac{d}{dt} \oa^t = \eta. $$ \\[-0.2cm]
Now we want to choose a 1-form $\si$ such that $\eta = d\si$. Define 
\renewcommand{\theequation}{4.8}
\begin{align}
\si(\ga,r) = - \int_{0}^{r} \eta_1 (\ga,s) ds.
\end{align}
By the fact that $|\eta|_{g_V} = |\oa - \oa_V|_{g_V} = O(r^{\nu})$ as $r \rightarrow 0$, we have $|\eta_0 |_{g_V} = |\eta_1 |_{g_V} = O(r^{\nu})$ as $r \rightarrow 0$, and using $\nu > 0$, we deduce that $\si$ is well-defined and it satisfies $|\si|_{g_V} = O(r^{\nu + 1})$ as $r \rightarrow 0$. Since $d\eta = 0$, we have $d_{\Ga} \eta_0 = 0$ and \\[-0.3cm]
\renewcommand{\theequation}{4.9}
\begin{align}
\frac{\pa \eta_0}{\pa r} + d_{\Ga} \eta_1 = 0, \\[-0.5cm] \notag
\end{align}
where $d_{\Ga}$ denotes the exterior differentiation in the $\ga$ direction. Therefore, \\[-0.3cm]
\begin{align*}
d\si(\ga,r) &= - \int_{0}^{r} d_{\Ga} (\eta_1 (\ga,s)) ds - dr \we \frac{\pa}{\pa r} \Big( \int_{0}^{r} \eta_1 (\ga,s) ds \Big) \\
&= \int_{0}^{r} \frac{\pa \eta_0}{\pa s} (\ga,s) ds - dr \we \eta_1 (\ga,r) \quad \textnormal{by (4.9).} \\[-0.5cm] \notag 
\end{align*}
For each fixed $r$, $\eta_{0}^{r} (\ga) := \eta_{0} (\ga, r)$ is a 2-form on $\Ga \cong \Ga \times \{1\}$, and so $|\eta_{0}^{r} (\ga)|_{g_V|_{\Ga \times \{1\}}} = r^2 |\eta_{0} (\ga, r)|_{g_V} = O(r^{\nu + 2 })$ since $|\eta_0|_{g_V} = O(r^{\nu})$. It follows that $|\eta_{0}^{r} (\ga)|_{g_V|_{\Ga \times \{1\}}} \rightarrow 0$ as $r \rightarrow 0$ and hence $\eta_{0} (\ga, r) \rightarrow 0$ as $r \rightarrow 0$. We then deduce that $d\si(\ga,r) = \eta_0 (\ga,r) + \eta_1 (\ga,r) \we dr = \eta(\ga,r)$. Therefore, we obtain a 1-form $\si$ on $\Ga \times (0, \ep)$ such that 
\begin{align*}
\frac{d}{dt} \oa^t = d\si. \\[-0.5cm]
\end{align*}
Also, $|\nab^k \si|_{g_V} = |\nab^{k-1} d\si|_{g_V} = |\nab^{k-1} \eta|_{g_V} = O(r^{\nu -k+1})$ as $r \rightarrow 0$ and for all $k \geq 1$. Together with $|\si|_{g_V} = O(r^{\nu + 1})$, we get a 1-form $\si$ satisfying 
\renewcommand{\theequation}{4.10}
\begin{align}
|\nab^k \si|_{g_V} = O(r^{\nu -k+1}) \quad \textnormal{as $r \rightarrow 0$ and for all $k \geq 0$.} \\ \notag
\end{align}
\par Now define a family of vector fields $X_t$ via \\[-0.5cm]
\renewcommand{\theequation}{4.11}
\begin{align}
\si + \iota(X_t) \oa^t = 0. \\[-0.5cm] \notag
\end{align}
The flow of this family of vector fields yields a family of diffeomorphisms $\psi_t$ on $V$ such that $\psi^{*}_{t} (\oa^t) = \oa^0$. 
In particular, we have constructed $\psi_1 : \Ga \times (0, \ep') \longrightarrow \Ga \times (0, \ep)$ for some $\ep'\in (0, \ep)$ which is a diffeomorphism with its image 
satisfying
$$\psi^{*}_{1} (\oa) = \psi^{*}_{1} (\oa^1) = \oa^0 = \oa_V. $$
Write $\Phi'= \Phi \circ \psi_1$ and $S' =  \Phi \circ \psi_1 (\Ga \times (0, \ep'))$, then $\Phi' : \Ga \times (0, \ep') \longrightarrow S'$ is a diffeomorphism such that
$$ (\Phi')^* (\oa_0) = \psi_{1}^{*} ( \Phi^* (\oa_0) ) = \psi_{1}^{*} (\oa) = \oa_V, $$
as required.  \\

\par From (4.11) we have $|\nab^k X_t |_{g_V} = O(r^{\nu -k+1})$ for all $k \geq 0$. Roughly speaking, $\psi_t$ = Id + $\int_{0}^{t} X_s ds$ to the first order, and hence $| \nab^k \psi_t |_{g_V} = O(r^{\nu - k +1})$ for all $k \geq 0$. It doesn't exactly make sense as $\psi_t$ and Id map to different points on $V$. But we could express them in terms of local coordinates $(x_1, \dots , x_{2m-1}, r)$ on $\Ga \times (0,\ep')$. Let $\psi^{j}_{t} (x_1, \dots , x_{2m-1}, r)$ be the $j$-th component function of $\psi_t$ for $j = 1, \dots, 2m$. Then $\pa^k (\psi^{j}_{t} (x_1, \dots , x_{2m-1}, r) - x_j ) = O(r^{\nu -k+ 1})$ for $j = 1, \dots, 2m-1$ and for all $k \geq 0$, and $\pa^k (\psi^{2m}_{t} (x_1, \dots , x_{2m-1}, r) - r) = O(r^{\nu -k+ 1})$ for all $k \geq 0$ where $\pa$ denotes the usual partial differentiation at the point $(x_1, \dots , x_{2m-1}, r)$. It follows that $\pa^k (\psi_t -$Id)$(x_1, \dots , x_{2m-1}, r) = O(r^{\nu - k + 1 })$ for all $k \geq 0$. Consequently we have 
$$|\pa^k (\psi_t - \textnormal{Id})^*(\Phi^*(\Oa_0))|_{g_V} = O(r^{\nu - k })$$ 
at the point $(x_1, \dots , x_{2m-1}, r)$. As a result, we have at each point on $\Ga \times (0, \ep')$
\begin{align*}
\left| \nab^k ( (\Phi')^* (\Oa_0) - \Oa_V ) \right|_{g_V} &= \left| \nab^k ( \psi_{1}^{*} ( \Phi^* (\Oa_0) ) - \Oa_V ) \right|_{g_V} \\
&\leq \left| \pa^k ( \psi_{1}^{*} ( \Phi^* (\Oa_0) ) - \Phi^* (\Oa_0) ) \right|_{g_V} \\
&\qquad + \left| \pa^k ( \Phi^* (\Oa_0) - \Oa_V ) \right|_{g_V} \\
&= O(r^{\nu - k}) + O(r^{\nu - k}) = O(r^{\nu - k})
\end{align*}
for all $k \geq 0$. This completes the proof. \hfill $\Box$ \\[1cm]

\begin{center}
\bf{4.3. \ Asymptotically Conical \CY \bf{\emph{m}}-folds} \\[0.3cm]
\end{center}

\normalsize

\par In the last part we study Asymptotically Conical (AC) \CY $m$-folds. We shall provide some examples and give an analogue of Theorem 4.8 for AC \CY $m$-folds. \\

\renewcommand{\thedefinition}{4.9}
\begin{definition}
\quad \textnormal{Let $(Y, J_Y, \oa_Y, \Oa_Y)$ be a complete, nonsingular \CY $m$-fold. Then $Y$ is an \emph{Asymptotically Conical (AC) \CY $m$-fold} with rate $\lda < 0$ modelled on a \CY cone $(V, J_V, \oa_V, \Oa_V)$ if there exist a compact subset $K \subset Y$, and a diffeomorphism $\Ups: \Ga \times (R, \infty) \longrightarrow Y \setminus K$ for some $R > 0$ such that
\begin{align}
|\nab^k ( \Ups^* ( \oa_Y ) - \oa_V ) |_{g_V} &= O(r^{\lda - k}), \ \textnormal{and} \tag{4.12} \\ \tag{4.13}
|\nab^k ( \Ups^* ( \Oa_Y ) - \Oa_V ) |_{g_V} &= O(r^{\lda - k}) \ \ \textnormal{as } r \rightarrow \infty \ \textnormal{and for all } k \geq 0. 
\end{align}
Here $\nab$ and $| \cdot |$ are computed using the cone metric $g_V$.} \\
\end{definition}
 
\par Similar asymptotic conditions on $g_Y$ and $J_Y$ can be deduced from (4.12) and (4.13). Unlike the conical singularity case, $[ \Ups^* ( \oa_Y ) - \oa_V ]$ and $[ \Ups^* ( \Oa_Y ) - \Oa_V ]$ need not be zero cohomology classes. Here are some conditions: \\

\renewcommand{\thelem}{4.10}
\begin{lem}
\quad Let $(Y, J_Y, \oa_Y, \Oa_Y)$ be an AC \CY $m$-fold with rate $\lda < 0$ modelled on a \CY cone $(V, J_V, \oa_V, \Oa_V)$. If $\lda < -2$ or $H^2 (\Ga, \BR) = 0$, then $[ \Ups^* ( \oa_Y ) - \oa_V ] = 0$. If $\lda < -m$ or $H^m (\Ga, \BC) = 0$, then $[ \Ups^* ( \Oa_Y ) - \Oa_V ] = 0$. \\
\end{lem}

\par The proof of it is similar to that of Lemma 4.7, except we now have $O(r^{\lda +2})$ for the integral of the difference of symplectic forms and $O(r^{\lda + m})$ for the holomorphic $(m,0)$-forms. Hence if $\lda < -2$, the integral approaches 0 as $r \rightarrow \infty$, which implies $[ \Ups^* ( \oa_Y ) - \oa_V ] = 0$. The same argument applies to the holomorphic $(m,0)$-forms. \\

\par We shall normally consider the case $\lda < -2$, so that $\Ups^* ( \oa_Y ) - \oa_V$ is always exact. Moreover, when $\lda < -2$, the proof for the analogue of Darboux Theorem works similarly as that for the conical singularities case. It is not clear whether the theorem holds for $\lda \geq -2$ and $[ \Ups^* ( \oa_Y ) - \oa_V ] = 0$ or not. \\

\renewcommand{\thethm}{4.11}
\begin{thm}
\quad Let $(Y, J_Y, \oa_Y, \Oa_Y)$ be an AC \CY $m$-fold with rate $\lda < -2$ modelled on a \CY cone $(V, J_V , \oa_V , \Oa_V)$. Then there exist a $R' > 0$ and a diffeomorphism $\Ups' : \Ga \times (R', \infty) \longrightarrow Y \setminus K$ such that $(\Ups')^* (\oa_Y) = \oa_ V$ and $|\nab^k ( (\Ups')^* ( \Oa_Y ) - \Oa_V ) |_{g_V} = O(r^{\lda - k})$ for all $k \geq 0$.  \\
\end{thm}

\par One can prove it in the same way as the proof of Theorem 4.8. The condition on the rate $\lda$ is essential for this proof to work. Since $|\eta_1|_{g_V} = O(r^{\lda})$, we need $\lda < -2$ to construct the 1-form $\si$. Moreover, in proving Theorem 4.8, we encountered the norm of $\eta_{0}^{r}$: $|\eta_{0}^{r} (\ga)|_{g_V|_{\Ga \times \{1\}}} = r^2 |\eta_{0} (\ga, r)|_{g_V}$, which is equal to $O(r^{\lda + 2 })$ in this case. Therefore we need $\lda < -2$ in order to have $\eta_0 \rightarrow 0$ as $r \rightarrow \infty$. \\[1cm]
\textbf{Example 4.12.}\quad Let $G$ be a finite subgroup of SU($m$) acting freely on $\BC^m \setminus \{0\}$, and $(X, \pi)$ a crepant resolution of the \CY cone $V = \BC^m / G$ given in Example 4.3. Then in each \ka class of ALE \ka metrics on $X$ there is a unique ALE Ricci-flat \ka metric (see Joyce \cite{Joyce1}, Chapter 8) and $X$ is then an AC \CY $m$-fold asymptotic to the cone $\BC^m / G$. In this case, it follows from \cite[Thm. 8.2.3]{Joyce1} that the rate $\lda$ is $-2m$. \\[-0.2cm]
\par If we take $G = \mathbb{Z}_m$ acting on $\BC^m$ as in Example 4.3, then a crepant resolution is given by the blow-up of $\BC^m / \mathbb{Z}_m$ at 0, which is also the total space of the canonical line bundle over $\BC \BP^{m-1}$. An explicit ALE Ricci-flat \ka metric is given in \cite[p.284-5]{Calabi} and also in \cite[Example 8.2.5]{Joyce1}. \\[1cm]
\textbf{Example 4.13.}\quad Consider the \CY cone $V = \{ z_{1}^2 + z_{2}^2 + z_{3}^2 +z_{4}^2 = 0 \}$ described in Example 4.4. One can construct two kinds of AC \CY manifolds. The first one is called the \emph{small resolution} of $V$, given by 
$$\widetilde{V} = \big\{ ( (z_1, \dots, z_4), [w_1, w_2]) \in \BC^4 \times \BC \BP^1 : z_1 w_2 = z_4 w_1, \ z_3 w_2 = z_2 w_1 \big\}. $$
It is essentially isomorphic to the normal bundle $\mathcal{O}(-1) \oplus \mathcal{O}(-1)$ over $\BC \BP^1$ with fibre $\BC^2$, and is also isomorphic to $V$ away from the origin where it is replaced by the whole $\BC \BP^1$. Note that one can obtain a second small resolution by swapping $z_3$ and $z_4$ in $\widetilde{V}$. Candelas et al \cite[p.258]{Can} constructed \CY metrics on $\widetilde{V}$, and it is an AC \CY 3-fold with rate $-2$. \\
\par The other is known as the \emph{deformation} or \emph{smoothing}, where $V$ is deformed to $Q_{\ep} = \{ z_{1}^2 + z_{2}^2 + z_{3}^2 +z_{4}^2 = \ep^2\}$ with $\ep$ a nonzero constant. This has the effect of replacing the node by an $S^3$. An important fact is that there is a symplectomorphism which identifies the standard symplectic form on $\BC^4$ restricted to $Q_{\ep}$ and the canonical symplectic form on the cotangent bundle $T^* S^3$ of $S^3$. Stenzel \cite[p.161]{Ste} constructed a \CY metric on $T^* S^3$, which makes $Q_{\ep}$ into an AC \CY 3-fold with rate $\lda = -3$.   \\[1cm]
\textbf{Example 4.14.}\quad Calabi \cite[p.284-5]{Calabi} constructed a 1-parameter family of AC \CY metrics on the canonical bundle $K_S$ of any K\"{a}hler-Einstein $(m-1)$-fold $S$ with positive scalar curvature, so that $K_S$ is an AC \CY $m$-fold modelled on the \CY cone $K_S \setminus S$ with rate $\lda = -2m$. \\
\par For the case in Example 4.5 (i), the $\mathcal{O}(-2,-2)$-bundle is AC \CY asymptotic to the cone $\mathcal{O}(-2,-2) \setminus (\BC \BP^1 \times \BC \BP^1)$ with rate $-6$. \\
\par Note that if we take $S = \BC \BP^{m-1}$, then we recover the case in Example 4.12 with $G = \mathbb{Z}_m$. \\[1cm]



\begin{center}
\large \bf{5. \ Calabi--Yau desingularizations} \\[0.5cm]
\end{center}





\normalsize

\par This section studies desingularizations of a compact \CY 3-fold $M_0$ with a conical singularity using an AC \CY 3-fold $Y$ with rate $\lda$. We shall only treat the simplest case here, in which $\lda < -3$ so that $\Ups_{t}^{*} (t^3 \Oa_Y) - \Oa_V$ is exact by Lemma 4.10. We explicitly construct a 1-parameter family of diffeomorphic, nonsingular compact 6-folds $M_t$ for small $t$ in $\S$5.1.  
Then in $\S$5.2 we construct a real closed 2-form $\oa_t$ and a complex closed 3-form $\Oa_t$ on $M_t$ and show that they give nearly \CY structures on $M_t$ for small enough $t$. Section 5.3 contains the main result of this paper, in which we show that the nearly \CY structure $(\oa_t, \Oa_t)$ on $M_t$ can be deformed to a genuine \CY structure $(\tilde{\oa}_t, \tilde{\Oa}_t)$ for small $t$ by applying Theorem 3.10. Finally in $\S$5.4, we apply our result to some examples studied before. We shall also discuss the case when $\lda = -3$.  \\

\begin{center}
\bf{5.1. \ Construction of $M_t$} \\[0.3cm]
\end{center}

\par Let $(M_0 , J_0 , \oa_0 , \Oa_0 )$ be a compact \CY 3-fold with a conical singularity at $x$ with rate $\nu$ modelled on a \CY cone $(V, J_V , \oa_V, \Oa_V)$. By Theorem 4.8, there exists an $\ep > 0$, a small open neighbourhood $S$ of $x$ in $M_0$ and a diffeomorphism $\Phi : \Ga \times (0, \ep) \longrightarrow S \setminus \{x\}$ such that $\Phi^* (\oa_0 ) = \oa_V$. \\

\par Let $(Y, J_Y , \oa_Y, \Oa_Y)$ be an AC \CY 3-fold with rate $\lda < -3$ modelled on the same \CY cone $V$. Theorem 4.11 shows that there is a diffeomorphism $ \Ups : \Ga \times (R, \infty) \longrightarrow Y \setminus\!K$ for some $R > 0$ such that 
$$\Ups^* (\oa_Y) = \oa_V \quad \textnormal{and} \quad |\nab^k (\Ups^* (\Oa_Y) - \Oa_V ) |_{g_V} = O(r^{\lda - k})$$ 
as $r \rightarrow \infty$ for all $k \geq 0$. We then apply a homothety to $Y$ such that 
$$ (Y, J_Y , \oa_Y, \Oa_Y) \longmapsto (Y, J_Y , t^2 \oa_Y,  t^3 \Oa_Y). $$
Then $(Y, J_Y , t^2 \oa_Y,  t^3 \Oa_Y)$ is also an AC \CY 3-fold, with the diffeomorphism $\Ups_t : \Ga \times (tR, \infty) \longrightarrow Y \setminus\!K$ given by  
$$ \Ups_t (\ga, r) = \Ups (\ga, t^{-1}r).$$
Our goal is to desingularize $(M_0 , J_0 , \oa_0 , \Oa_0 )$ by gluing $(Y, J_Y , t^2 \oa_Y,  t^3 \Oa_Y)$ in at $x$ to produce a family of compact nonsingular \CY 3-folds. \\

\par Fix $\al \in (0,1)$ and let $t > 0$ be small enough that $tR < t^{\al} < 2t^{\al} < \ep$. Define
\begin{align*}
P_t &= K \cup \Ups_t ( \Ga \times (tR, 2t^{\al}) ) \subset Y \quad \textnormal{and} \\
Q_t &= (M_0 \setminus S) \cup \Phi(\Ga \times (t^{\al}, \ep) ) \subset M_0 .
\end{align*} 
The diffeomorphism $\Phi \circ \Ups_{t}^{-1}$ identifies $\Ups_t ( \Ga \times (t^{\al}, 2t^{\al}) ) \subset P_t$ and $\Phi(\Ga \times (t^{\al}, 2t^{\al}) ) \subset Q_t$, and we define the intersection $P_t \cap Q_t$ to be the region $\Ups_t ( \Ga \times (t^{\al}, 2t^{\al}) ) \cong \Phi(\Ga \times (t^{\al}, 2t^{\al}) ) \cong \Ga \times (t^{\al}, 2t^{\al})$. Define $M_t$ to be the quotient space of the union $P_t \cup Q_t$ under the equivalence relation identifying the two annuli $\Ups_t ( \Ga \times (t^{\al}, 2t^{\al}) )$ and $\Phi(\Ga \times (t^{\al}, 2t^{\al}) )$. Then $M_t$ is a smooth nonsingular compact 6-fold for each $t$. \\[1cm]

\begin{center}
\bf{5.2. \ Nearly \CY structures $(\oa_t, \Oa_t)$: the case $\lda < -3$ }\\[0.3cm] 
\end{center}

\par In this section we construct on $M_t$ a real closed 2-form $\oa_t$ and a complex closed 3-form $\Oa_t$, and show  they together give nearly \CY structures on $M_t$ for small enough $t$. Define
\begin{center}
$\oa_t$ = 
$\begin{cases}$
$\oa_0$ \ \ on\ \ $Q_t$, $\\$
$t^2 \oa_Y$\ \ on\ \ $P_t$.
$\end{cases}$
\end{center}

This is well-defined as $\Phi^* (\oa_0) = \oa_V = \Ups_{t}^{*} (t^2 \oa_Y)$ on the intersection $P_t \cap Q_t$ by Theorem 4.8 and 4.11. Thus $\oa_t$ gives a symplectic form on $M_t$. \\

\par Let $F: \BR \longrightarrow [0,1]$ be a smooth, increasing function with $F(s) = 0$ for $s \leq 1$ and $F(s) = 1$ for $s \geq 2$. Then for $r \in (tR, \ep)$, $F(t^{-\al}r) = 0$ for $tR < r \leq t^{\al}$ and $F(t^{-\al}r) = 1$ for $2t^{\al} \leq r < \ep$. We now define a complex 3-form on $M_t$. From (4.3), we have $|\Phi^* (\Oa_0) - \Oa_V |_{g_V} = O(r^{\nu})$. As $\nu > 0$, it follows that $\Phi^* (\Oa_0) - \Oa_V$ is exact, and we can write 
\renewcommand{\theequation}{5.1}
\begin{align}
\Phi^* (\Oa_0) = \Oa_V + dA
\end{align}
for some complex 2-form $A(\ga, r)$ on $\Ga \times (0, \ep)$ satisfying 
\renewcommand{\theequation}{5.2}
\begin{align}
| \nab^k A(\ga, r)|_{g_V} = O(r^{\nu + 1 -k}) \quad \textnormal{as } r \rightarrow 0 \  \textnormal{ for all } k \geq 0. 
\end{align}
The case $k = 0$ follows by defining $A$ by integration as in Theorem 4.8. Similarly, as we have assumed $\lda < -3$ to simplify the problem, the 3-form $\Ups^* (\Oa_Y) - \Oa_V$ is exact by Lemma 4.10 and we can write
$$\Ups^* (\Oa_Y) = \Oa_V + dB $$
for some complex 2-form $B(\ga, r)$ on $\Ga \times (R, \infty)$ satisfying 
$$| \nab^k B(\ga, r)|_{g_V} = O(r^{\lda + 1 -k}) \quad \textnormal{as } r \rightarrow \infty \textnormal{ and for all } k \geq 0. $$ 
Then we apply a homothety to $Y$ and rescale the forms to get $B(\ga, t^{-1}r)$ on $\Ga \times (tR, \infty)$ such that 
\renewcommand{\theequation}{5.3}
\begin{align}
\Ups_{t}^{*} (t^3 \Oa_Y) = \Oa_V + t^3 dB(\ga, t^{-1}r)
\end{align}
and
\renewcommand{\theequation}{5.4}
\begin{align}
| \nab^k B(\ga, t^{-1} r)|_{g_V} = O(t^{-\lda-3} r^{\lda + 1 -k}) \quad \textnormal{for } r > tR \textnormal{ and for all } k \geq 0. 
\end{align}
Define a smooth, complex closed 3-form $\Oa_t$ on $M_t$ by \\[-0.3cm]
\renewcommand{\theequation}{5.5}
\begin{align}
\Oa_t = 
\begin{cases}
\Oa_0$ \ \ on\ \ $Q_t \setminus (P_t \cap Q_t), \\
\Oa_V + d\big[ F(t^{-\al}r) A(\ga, r) + t^3 (1 - F(t^{-\al}r)) B(\ga, t^{-1}r)\big]$ \ on\ $P_t \cap Q_t, \\ 
t^3 \Oa_Y$\ \ on\ \ $P_t \setminus (P_t \cap Q_t).   
\end{cases} \\[-0.3cm]  \notag
\end{align} 
Note that when $2t^{\al} \leq r < \ep$ we have $F(t^{-\al}r) = 1$ so that $\Oa_t = \Phi^* (\Oa_0)$ by (5.1), and when $tR < r \leq t^{\al}$ we have $F(t^{-\al}r) = 0$, so that $\Oa_t = \Ups_{t}^{*} (t^3 \Oa_Y)$ by (5.3). Therefore, $\Oa_t$ interpolates between $\Phi^* (\Oa_0)$ near $r = \ep$ and $\Ups_{t}^{*} (t^3 \Oa_Y)$ near $r = tR$. \\

\par Recall that if we get a real closed 2-form $\oa$ and a complex 3-form $\Oa$ which are sufficiently close to the \ka form $\tilde{\oa}$ and the holomorphic volume form $\tilde{\Oa}$ of a \CY structure respectively, then Proposition 3.2 tells us that $(\oa, \Oa)$ gives a nearly \CY structure on the manifold. Making use of this idea, we have \\[-0.3cm]

\renewcommand{\theprop}{5.1} 
\begin{prop}
\quad Let $M_t$, $\oa_t$ and $\Oa_t$ be defined as above. Then $(\oa_t, \Oa_t)$ gives a nearly \CY structure on $M_t$ for sufficiently small t. \\[-0.3cm]
\end{prop}
\pf \ \, We only have to prove the statement on $P_t \cap Q_t$, as $(M_t, \oa_t, \Oa_t)$ is \CY on $P_t \setminus P_t \cap Q_t$\, and on \,$Q_t \setminus P_t \cap Q_t$, and hence is nearly Calabi-Yau. We prove it by applying Proposition 3.2, that is, we show on $P_t \cap Q_t$ that $(\oa_t, \Oa_t)$ is sufficiently close to the genuine \CY structure $(\oa_V, \Oa_V)$ coming from the \CY cone $V$ for small $t$. We choose to compare with $(\oa_V, \Oa_V)$ rather than either of the \CY structures $(\oa_0, \Oa_0)$ and $(t^2\oa_Y, t^3\Oa_Y)$ on $P_t \cap Q_t$ since we have already got bounds on norms for various forms w.r.t. the cone metric $g_V$. Now $\oa_t = \oa_V$\, on \,$P_t \cap Q_t$, while 
$$ \Oa_t - \Oa_V \,=\, d\big[ F(t^{-\al}r) A(\ga, r) + t^3 (1 - F(t^{-\al}r)) B(\ga, t^{-1}r)\big] \ \textnormal{ on }\ P_t \cap Q_t  $$
by (5.5). Calculation shows that 
\renewcommand{\theequation}{5.6}
\begin{align}
|(\Oa_t - \Oa_V )(\ga, r)|_{g_V} = O(t^{-\lda(1 - \al)}) + O(t^{\al \nu}) \quad \textnormal{for } r \in (t^{\al}, 2t^{\al}),
\end{align}
and hence $|\Oa_t - \Oa_V |_{g_V} \leq C_0 t^{\ga}$ where $C_0 > 0$ is some constant and $\ga$ = min($-\lda(1 - \al)$, $\al \nu$). Hence Proposition 3.2 applies with $\ep = C_0 t^{\ga}$ if $t$ is small enough such that $C_0 t^{\ga} \leq \ep_1$, and so $(\oa_t, \Oa_t)$ gives a nearly \CY structure on $P_t \cap Q_t$. This completes the proof. \hfill $\Box$ \\

\par Therefore we can associate an almost complex structure $J_t$ and a real 3-form $\ta_{2,t}'$ such that $\Oa_t' :=$ Re($\Oa_t$) + $i \ta_{2,t}'$ is a (3,0)-form w.r.t. $J_t$. Moreover, we have the 2-form $\oa_t'$, which is the rescaled (1,1)-part of $\oa_t$ w.r.t. $J_t$, and the associated metric $g_t$ on $M_t$. Following similar arguments to Proposition 3.2, we conclude that $|g_t - g_V|_{g_V} = O(t^{-\lda(1 - \al)}) + O(t^{\al \nu}) = |g^{-1}_{t} - g^{-1}_{V}|_{g_V}$ . \\


\begin{center}
\bf{5.3. \ The main result} \\[0.3cm]
\end{center}

\par We are now ready to prove our main result on desingularization of compact \CY 3-folds $M_0$ with a conical singularity in the simplest case when $\lda < -3$. The proof of it uses Theorem 3.10. \\


\renewcommand{\thethm}{5.2} 
\begin{thm}
\quad Suppose $(M_0, J_0, \oa_0, \Oa_0)$ is a compact \CY 3-fold with a conical singularity at $x$ with rate $\nu > 0$ modelled on a \CY cone $(V, J_V , \oa_V, \Oa_V)$. Let $(Y, J_Y , \oa_Y, \Oa_Y)$ be an AC \CY 3-fold with rate $\lda < -3$ modelled on the same \CY cone $V$. Define a family $(M_t, \oa_t, \Oa_t)$ of nonsingular compact nearly \CY 3-folds, with \CY metrics $g_t$ as in $\S$5.1 and $\S$5.2. 
\par Then $M_t$ admits a \CY structure $(\tilde{J}_t, \tilde{\oa}_t, \tilde{\Oa}_t)$ such that $\| \tilde{\oa}_t - \oa_t \|_{C^0} \leq K t^{\kappa}$ and\, $\| \tilde{\Oa}_t - \Oa_t \|_{C^0} \leq K t^{\kappa}$ for some $\kappa, K > 0$ and for sufficiently small $t$. The cohomology classes satisfy $[\textnormal{Re}(\Oa_t)]$ = $[\textnormal{Re}(\tilde{\Oa}_t)] \in H^3 (M_t, \BR)$\, and\, $[\oa_t]$ = $c_t\,[\tilde{\oa}_t] \in H^2 (M_t, \BR)$ for some $c_t > 0$. Here all norms are computed with respect to $g_t$. \\
\end{thm}
\pf \, First we estimate the norms of $\oa_t - \oa_t'$ and $\textnormal{Im}(\Oa_t) - \ta_{2,t}' = \textnormal{Im}(\Oa_t) - \textnormal{Im}(\Oa_t')$ on $P_t \cap Q_t$, as in part (i) of Theorem 3.10. Since $\oa_t'$ depends on Re($\Oa_t$) and $\oa_t$ (= $\oa_V$ on $P_t \cap Q_t$) on Re($\Oa_V$), it follows that 
\begin{align*}
|\oa_t - \oa_t'|_{g_t} &\leq C_1 |\oa_t - \oa_t'|_{g_V} \leq C_2 |\textnormal{Re}(\Oa_t) - \textnormal{Re}(\Oa_V) |_{g_V} \leq C_2 | \Oa_t - \Oa_V |_{g_V} \\
&= O(t^{-\lda(1 - \al)}) + O(t^{\al \nu})
\end{align*}
for some constants $C_1, C_2 > 0$ and hence 
\renewcommand{\theequation}{5.7}
\begin{align}
\| \oa_t - \oa_t'\|_{C^0} = O(t^{-\lda(1 - \al)}) + O(t^{\al \nu}). 
\end{align}
From the fact that vol($P_t \cap Q_t$) = $O(t^{6\al})$, we have 
\renewcommand{\theequation}{5.8}
\begin{align}
\| \oa_t - \oa_t'\|_{L^2} = O(t^{\frac{6\al}{2}}) \cdot \| \oa_t - \oa_t'\|_{C^0} = O(t^{3\al-\lda(1 - \al)}) + O(t^{3\al+\al \nu}).
\end{align}
Furthermore, 
\begin{align*}
|\textnormal{Im}(\Oa_t) - \textnormal{Im}(\Oa_t')|_{g_t} &\leq C_3 \,|\textnormal{Im}(\Oa_t) - \textnormal{Im}(\Oa_t')|_{g_V} \\
&\leq C_3 \,|\textnormal{Im}(\Oa_t) - \textnormal{Im}(\Oa_V)|_{g_V} + C_3 \,|\textnormal{Im}(\Oa_V) - \textnormal{Im}(\Oa_t')|_{g_V} \\
&\leq C_3 \,|\Oa_t - \Oa_V|_{g_V} + C_4 \,|\textnormal{Re}(\Oa_V) - \textnormal{Re}(\Oa_t) |_{g_V} \\
&= O(t^{-\lda(1 - \al)}) + O(t^{\al \nu})
\end{align*}
for some constants $C_3, C_4 > 0$, as Im($\Oa_V$) is determined by Re($\Oa_V$) and Im($\Oa_t)'$ by Re($\Oa_t$). Therefore,
\renewcommand{\theequation}{5.9}
\begin{align}
\|\textnormal{Im}(\Oa_t) - \textnormal{Im}(\Oa_t')\|_{C^0} = O(t^{-\lda(1 - \al)}) + O(t^{\al \nu})
\end{align}
and
\renewcommand{\theequation}{5.10}
\begin{align}
\|\textnormal{Im}(\Oa_t) - \textnormal{Im}(\Oa_t')\|_{L^2} = O(t^{3\al-\lda(1 - \al)}) + O(t^{3\al+\al \nu}). \\ \notag
\end{align}
\par It can be deduced from (5.5) and (5.6) that $| \nab^{g_t} (\Oa_t - \Oa_V) |_{g_t} = O(t^{-\lda(1 - \al) - \al}) + O(t^{\al \nu - \al})$ and $| (\nab^{g_t})^2 (\Oa_t - \Oa_V) |_{g_t} = O(t^{-\lda(1 - \al) - 2\al}) + O(t^{\al \nu - 2\al})$, which imply the equalities
\renewcommand{\theequation}{5.11}
\begin{align}
\| \nab^{g_t} (\oa_t - \oa_t') \|_{C^0} = O(t^{-\lda(1 - \al) - \al}) + O(t^{\al \nu - \al})
\end{align}
and 
\renewcommand{\theequation}{5.12}
\begin{align}
\| (\nab^{g_t})^2 (\oa_t - \oa_t') \|_{C^0} = O(t^{-\lda(1 - \al) - 2\al}) + O(t^{\al \nu - 2\al}).
\end{align}
Then the $L^{12}$-norm satisfies
\renewcommand{\theequation}{5.13}
\begin{align}
\| \nab^{g_t} (\oa_t - \oa_t') \|_{L^{12}} &= O(t^{\frac{6\al}{12}}) \cdot \| \nab^{g_t} (\oa_t - \oa_t') \|_{C^0} \notag \\
&= O(t^{-\frac{\al}{2} - \lda(1 - \al)}) + O(t^{-\frac{\al}{2} + \al \nu }).
\end{align}
Finally, we estimate the $L^{12}$-norms of $\nab^{g_t} \oa_t$ and $\nab^{g_t} \textnormal{Re}(\Oa_t)$. Note that
$$|\nab^{g_t} \oa_t|_{g_t} \leq C_5 \,|(\nab^{g_t} - \nab^{g_V}) \oa_t \,|_{g_V} + C_5 \,|\nab^{g_V} \oa_t|_{g_V} = C_5 \,|(\nab^{g_t} - \nab^{g_V}) \oa_t \,|_{g_V}$$
for some constant $C_5 > 0$, as $\oa_t = \oa_V$ on $P_t \cap Q_t$ and $\nab^{g_V} \oa_V = 0$. Then
\begin{align*}
|\nab^{g_t} \oa_t|_{g_t} &\leq C_5 \,|(\nab^{g_t} - \nab^{g_V})|_{g_V} \cdot \,|\oa_t|_{g_V} \\
&\leq \frac{3}{2}\,C_5 \,|g^{-1}_{t}|_{g_V} \cdot |\nab^{g_V} g_t|_{g_V} \cdot \,|\oa_t|_{g_V} \quad \textnormal{by (3.21)} \\
&= C_6 |\nab^{g_V} (g_t - g_V) |_{g_V} \quad \textnormal{as } \nab^{g_V} g_V = 0.
\end{align*}
Here $C_6$ is an upper bound for $\frac{3}{2}\,C_5 \,|g^{-1}_{t}|_{g_V} \cdot \,|\oa_t|_{g_V}$ which is independent of $t$. It follows that 
\begin{align*}
|\nab^{g_t} \oa_t|_{g_t} = O(t^{-\lda(1 - \al) - \al}) + O(t^{\al \nu - \al}),
\end{align*}
and consequently
\renewcommand{\theequation}{5.14}
\begin{align}
\|\nab^{g_t} \oa_t \|_{L^{12}} = O(t^{-\frac{\al}{2} - \lda(1 - \al)}) + O(t^{-\frac{\al}{2} + \al \nu }).
\end{align}
A similar argument shows 
\renewcommand{\theequation}{5.15}
\begin{align}
\|\nab^{g_t}\textnormal{Re}(\Oa_t) \|_{L^{12}} = O(t^{-\frac{\al}{2} - \lda(1 - \al)}) + O(t^{-\frac{\al}{2} + \al \nu }).
\end{align}

\par Now for parts (i) to (iii) of Theorem 3.10 to hold, we need: \\

$\begin{cases}$
$-\lda(1 - \al) \geq \kappa, \quad \al \nu \geq \kappa \quad \textnormal{from (5.7) and (5.9),} \\
3\al-\lda(1 - \al) \geq 3+\kappa, \quad 3\al+\al \nu \geq 3+\kappa \quad \textnormal{from (5.8) and (5.10),} \\
-\frac{\al}{2} - \lda(1 - \al) \geq -\frac{1}{2} + \kappa, \quad -\frac{\al}{2} + \al \nu \geq -\frac{1}{2} + \kappa  \quad \textnormal{from (5.13), (5.14) and (5.15),}\\
-\lda(1 - \al) - \al \geq \kappa - 1, \quad \al \nu - \al \geq \kappa - 1 \quad \textnormal{from (5.11),} \\
-\lda(1 - \al) - 2\al \geq \kappa - 2, \quad \al \nu - 2\al \geq \kappa - 2 \quad \textnormal{from (5.12).} $
$\end{cases}$ \\[0.3cm]


\par Observe that the second set of inequalities imply all the others, as $\al \leq 1$. Therefore, calculations using these two inequalities show that there exist solutions $\al \in (0,1)$ and $\kappa > 0$ for any $\nu > 0$ and $\lda < -3$. For example, we could take
$$\al = \frac{1}{2}\left(\frac{6+\nu}{3+\nu}\right) \in (0,1) \quad \textnormal{and} \quad \kappa = \textnormal{min}\left((1-\al)(-3-\lda), \frac{\nu}{2}\right) > 0.$$ \\


\par For parts (iv) and (v) of Theorem 3.10, note that under the homothety $g_Y \mapsto t^2 g_Y$ on the AC \CY 3-fold $Y$ we have $\da(t^2 g_Y )= t\,\da(g_Y)$ and $\| R(t^2 g_Y) \|_{C^0} = t^{-2} \| R(g_Y) \|_{C^0}$. Moreover, the dominant contributions to $\da(g_t)$ and $\| R(g_t) \|_{C^0}$ for small $t$ come from $\da(g_Y)$ and $\| R(g_Y) \|_{C^0}$ which are proportional to $t$ and $t^{-2}$. Thus there exist constants $E_3, E_4 > 0$ such that (iv), (v) of Theorem 3.10 hold for sufficiently small $t$. Hence by Theorem 3.10, $M_t$ admits a \CY structure $(\tilde{J}_t, \tilde{\oa}_t, \tilde{\Oa}_t)$ such that $\| \tilde{\oa}_t - \oa_t \|_{C^0} \leq K t^{\kappa}$ and\, $\| \tilde{\Oa}_t - \Oa_t \|_{C^0} \leq K t^{\kappa}$ for some $\kappa, K > 0$ and for sufficiently small $t$.\\

\par Finally, the cohomology condition in Theorem 3.10 holds automatically here. This is because if $Y$ is an AC \CY 3-fold, then it can be shown that Hol($g_Y$) = SU(3) except for the trivial case where $Y = V = \BC^3$. Moreover, Hol($g_Y$) lies inside the holonomy group Hol($\tilde{g}_t$) of the \CY metric $\tilde{g}_t$ on $M_t$, and hence Hol($\tilde{g}_t$) must be the whole SU(3). The first cohomology group $H^1 (M_t, \BR)$ therefore vanishes for each sufficiently small $t$, and the theorem now follows from Theorem 3.10.  \hfill $\Box$ \\

\par Theorem 5.2 can be extended to the case when the \CY 3-fold $(M_0, J_0,$ $\oa_0, \Oa_0)$ has finitely many conical singularities at $x_1, \dots , x_n$ with rates $\nu_1, \dots ,\nu_n $ $> 0$ modelled on \CY cones $V_1, \dots , V_n$. Let $(Y_i, J_{Y_i}, \oa_{Y_i}, \Oa_{Y_i})$ be AC \CY 3-folds with rates $\lda_i < -3$ modelled on cones $V_i$ for $i = 1, \dots , n$. We then glue in at the singular points AC \CY 3-folds $(Y_i, J_{Y_i}, t^2 \oa_{Y_i}, t^3 \Oa_{Y_i})$ and produce a 1-parameter family of nonsingular \CY 3-folds. \\[0.3cm]

\begin{center}
\bf{5.4. \ Conclusions} \\[0.3cm]
\end{center}

\par We conclude by applying the above result to some examples given in $\S$3 and discussing briefly the case $\lda = -3$. First consider the situation in Example 4.12 and take $m = 3$. Then the crepant resolution $X$ of the \CY cone $V = \BC^3 / G$ is an AC \CY 3-fold with rate $\lda = -6$. Thus Theorem 5.2 applies and we can desingularize any compact \CY 3-fold $M_0$ with conical singularities modelled on $V$, or equivalently, any \CY 3-orbifold with isolated singularities, by the gluing process. In particular when $G = \mathbb{Z}_3$, a standard example of compact \CY 3-orbifold with isolated singularities is given in \cite[Example 6.6.3]{Joyce1}. Let $\Lambda = <1,\ e^{2\pi i/3}>$ be a lattice in $\BC$, and take $T^6 =  T^2 \times T^2 \times T^2$ where we regard each $T^2$ as the quotient $\BC / \Lambda$. The group $\mathbb{Z}_3$ acts on $T^6$ by multiplication by $e^{2\pi i/3}$ on each $T^2$-component. Then $M_0 = T^6 / \mathbb{Z}_3$ is a \CY 3-orbifold and it is not hard to see that $M_0$ has 27 isolated singular points modelled on $\BC^3 / \mathbb{Z}_3$. Thus by gluing in the total space $K_{\BC \BP^2} = \mathcal{O}(-3) \longrightarrow \BC \BP^2$ of the canonical bundle over $\BC \BP^2$ in each of the singular points, we obtain a desingularization of $M_0 = T^6 / \mathbb{Z}_3$. \\


\par Note that if we desingularize a \CY 3-orbifold with isolated singularities modelled on $\BC^3 / G$ by gluing, the \emph{Schlessinger Rigidity Theorem} \cite{Sch} tells us that $\BC^3 / G$ admits no nontrivial deformations, and hence what we will get will be a crepant resolution of the original orbifold. Potentially we will get a crepant resolution of a deformation of the orbifold, but it will be the crepant resolution of the original orbifold if [Re($\Oa_t$)] = [Re($\Oa$)] for all sufficiently small $t$, or equivalently if $[\Oa_Y] = 0$ in $H^3 (Y, \BC)$. \\

\par Now Yau's solution to the Calabi conjecture \cite{Yau} gives the existence of \CY metrics on the crepant resolution. However, it does not provide a way to write down the \CY metrics explicitly, and so in general we do not know much about what the \CY metrics are like. But in the orbifold case, our result tells a bit more by giving a quantitative description of these \CY metrics, showing that these metrics locally look like the metrics obtained by gluing the orbifold metrics and the ALE metrics on the crepant resolution of $\BC^3 / G$.  \\ 

\par Our result can also be applied to desingularize compact \CY 3-folds with conical singularities modelled on the \CY cone $\mathcal{O}(-2,-2) \setminus (\BC \BP^1 \times \BC \BP^1)$ by gluing in the AC \CY 3-fold $\mathcal{O}(-2,-2)$-bundle with rate $-6$. Thus we could resolve a kind of singularity which is not of orbifold type. \\

\par Finally we would like to discuss the case when $\lda = -3$. In Theorem 5.2 we desingularize a \CY 3-fold with a conical singularity using an AC \CY 3-fold with rate $\lda$ where we assumed that $\lda < -3$ and hence $[\Ups_{t}^{*} (t^3 \Oa_Y) - \Oa_V] = 0$ by Lemma 4.10. If we relax this to allow $\lda = -3$, then the cohomology class $[\Ups_{t}^{*} (t^3 \Oa_Y) - \Oa_V]$ may be nonzero in $H^3 (\Ga, \BC)$. But the cohomology class $[\Phi^* ( \Oa_0 ) - \Oa_V]$ is always zero by Lemma 4.7, so there can be \emph{topological obstructions} to defining a closed 3-form which interpolates between $\Phi^* ( \Oa_0 )$ and $\Ups_{t}^{*} (t^3 \Oa_Y)$. Thus allowing $\lda = -3$ introduces global problems for our gluing method. Here is a very short sketch of a method the author hopes to use to tackle the problem. We replace the holomorphic (3,0)-form $\Oa_0$ on $M_0$ by $\Oa_0 + t^3 \chi$, where $\chi$ is some closed and coclosed (2,1)-form with appropriate asymptotic behaviour, and $[\Phi^* (\chi)] = [\Ups^* (\Oa_Y)]$ on $H^3 (\Ga, \BC)$. Calculations by the author suggest that such $\chi$ exists and it cancels out the $O(t^3 r^{-3})$ terms such that Theorem 3.10 can handle the size of the error introduced. The advantage of extending our result to the case $\lda = -3$ is that it can be applicable to a larger class of AC \CY 3-folds such as the deformation $Q_{\ep}$ of the cone $V = \{ z_{1}^2 + z_{2}^2 + z_{3}^2 +z_{4}^2 = 0 \}$ in Example 4.13. Hence by extending our result, we could smooth ordinary double points by analytic rather than algebro-geometric methods. \\[1cm]  




\normalsize

\begin{center}

\end{center}






\end{document}